\documentclass[ijoc]{informs4} 

\usepackage{eqndefns-left} 
 
\RequirePackage{tgtermes}
\RequirePackage{newtxtext}
\RequirePackage{newtxmath}
\RequirePackage{bm}
\RequirePackage{endnotes}

\def\theARTICLETOP{}
\def\theLRHFirstLine{}
\def\theLRHSecondLine{}
\def\theRRHFirstLine{}
\def\theRRHSecondLine{}

\usepackage{multirow}
\usepackage{hyperref}

\definecolor{DarkNavy}{rgb}{0.0, 0.13, 0.33}
\definecolor{DarkMaroon}{rgb}{0.5, 0.0, 0.0}
\definecolor{DarkGreen}{rgb}{0.0, 0.3, 0.1}

\hypersetup{
    colorlinks=true,       
    linkcolor=DarkNavy,    
    citecolor=DarkMaroon,  
    urlcolor=DarkGreen,    
    linktoc=all            
}

\usepackage{algorithm}
\usepackage{algpseudocode}
\usepackage{xcolor}
\usepackage{etoolbox}

\AtBeginEnvironment{equation}{\small}
\AtBeginEnvironment{equation*}{\small}
\AtBeginEnvironment{align}{\small}
\AtBeginEnvironment{align*}{\small}

\newcommand{\xub}{\overline{x}}
\newcommand{\xlb}{\underline{x}}
\newcommand{\xhat}{\hat{x}}
\newcommand{\lamtil}{\tilde{\lambda}}
\newcommand{\lamub}{\overline{\lambda}}

\newcommand{\lamhat}{\hat{\lambda}}

\OneAndAHalfSpacedXII 

\usepackage{algorithm}
\usepackage{algpseudocode}
\usepackage{tikz}
\usepackage{natbib}
 \bibpunct[, ]{(}{)}{,}{a}{}{,}
 \def\bibfont{\small}

\EquationsNumberedThrough    
\TheoremsNumberedThrough     
\ECRepeatTheorems  
\MANUSCRIPTNO{IJOC- }

\begin{document}
\RUNAUTHOR{Wu, Luo, and Ouyang}

\RUNTITLE{Parameter-Free Constrained Variational Inequalities}
\TITLE{Parameter-Free First-Order Method for  Constrained Variational Inequalities}

\ARTICLEAUTHORS{
\AUTHOR{Yan Wu}
\AFF{Department of Industrial Engineering,
Clemson University, \EMAIL{yw8@clemson.edu }}

\AUTHOR{Qi Luo}
\AFF{Department of Business Analytics,
University of Iowa, \EMAIL{qi-luo-1@uiowa.edu}}

\AUTHOR{Yuyuan Ouyang}
\AFF{School of Mathematical and Statistical Sciences,
Clemson University, \EMAIL{ yuyuano@clemson.edu}}
} 

\ABSTRACT{
 Constrained variational inequalities (VIs) provide a general framework for modeling and solving complementarity and equilibrium problems with functional constraints. In this paper, we propose a parameter-free accelerated operator sliding (PFAOS) method for solving monotone  VIs  with convex smooth constraints. For a single constraint, PFAOS achieves an operator evaluation complexity of $\mathcal{O}(1/\varepsilon)$ and a constraint gradient evaluation complexity of $\mathcal{O}(1/\sqrt{\varepsilon})$, without requiring prior knowledge of the Lipschitz constant of the VI operator or the smoothness constant of the constraint function. We further extend PFAOS to multiple constraints by incorporating a primal-dual hybrid gradient update, yielding a multi-constrained PFAOS (MC-PFAOS) method that preserves the same evaluation complexities.  
 Finally, numerical experiments are conducted on large-scale games with coupling constraints to validate the efficiency and effectiveness of the proposed algorithms.
}

\FUNDING{This research was supported by the Air Force Office of Scientific Research grant FA9550-25-1-0278.}


\maketitle
\renewcommand{\theARTICLETOP}{}

\section{Introduction}
Many equilibrium models incorporate functional constraints to capture operational and physical restrictions on decentralized decisions.
For example, in transportation networks, Wardrop equilibrium route choices could yield network flows subject to traffic-control, network-wide emission, or battery-range constraints \citep{larsson1999side, zhang2019range}.
Similar structures arise in communication and computing systems, where users compete for bandwidth under quality-of-service constraints \citep{johari2004efficiency}. 
In supply-chain and oligopoly models, firms may compete in a Cournot equilibrium subject to limited supplier capacities or firm-level resource constraints \citep{cho2014capacity, caldentey2024multimarket}.
These problems can be formulated as constrained variational inequalities (CVIs), where equilibrium solutions are subject to many nonlinear functional constraints. These feasibility requirements pose substantial computational challenges beyond those of unconstrained equilibrium problems.

We consider the CVIs of finding
a point $x^*\in\widetilde X$ such that
\begin{equation}\label{eq:vi_strong_solution}
      \langle F(x^*),x^*-x\rangle \leq 0,
   \quad \forall x\in\widetilde X :=
   \{x\in X:G(x)\leq 0\}.
\end{equation}
Here, $X\subseteq\mathbb R^n$ is a compact convex set onto which projection is computationally inexpensive.
The operator $F:\mathbb R^n\to\mathbb R^n$ is monotone and Lipschitz continuous with respect to the Euclidean norm, that is, there exists $M>0$ such that
\begin{equation*}
       \langle F(u)-F(v),u-v\rangle \geq0,
    \qquad
    \|F(u)-F(v)\|
    \leq M\|u-v\|,
    \quad \forall u,v\in X.
\end{equation*}
The constraint mapping $G(x):=(g_1(x),\ldots,g_m(x))^\top$ consists of convex and differentiable functions $g_i:\mathbb R^n\to\mathbb R$, with $G(x)\leq 0$ interpreted componentwise.
$\|\cdot\|$ denotes the Euclidean norm for vectors and the induced operator norm for matrices
Its Jacobian  $J(x)$ has $\nabla g_i(x)^\top$ as its $i$th row. 
In addition, each $g_i$ has an $L_{(i)}$-Lipschitz continuous gradient on $X$, with $L_{(i)}>0$, so that
\begin{equation*}\label{eq:L_smooth}
    0\leq
    g_i(v)-g_i(u) -\langle\nabla g_i(u),v-u\rangle
    \leq \tfrac{L_{(i)}}{2}\|v-u\|^2, \quad \forall u,v\in X.
\end{equation*}
We denote the aggregate smoothness constant by
\(
L:=  (\sum_{i=1}^m L_{(i)}^2  )^{1/2}.
\)

A point $x^*$ satisfying \eqref{eq:vi_strong_solution} is called an exact \emph{strong} (Stampacchia) solution.
In contrast, an exact \emph{weak} (Minty) solution is a point $x^*\in\widetilde X$ satisfying
$\langle F(x),x^*-x\rangle\leq 0$ for all $x\in\widetilde X$.
The exact strong and weak solution sets coincide when $F$ is continuous and monotone. 
Our goal is to compute an $\varepsilon$-approximate weak solution, namely, a point  $\bar x \in X$ such that
\begin{equation} \label{eq:problem_of_interest}
     \sup\nolimits_{x\in  \widetilde X } \langle F(x),\bar x-x\rangle \leq \varepsilon,
     \text{ and }
     \|[G(\bar x)]_+\|\leq\varepsilon,
\end{equation}
where $[\cdot]_+$ denotes the componentwise maximum with zero. 

Let $ri(C)$ denote the relative interior of a set \(C\).
We further impose the following assumption such that the CVI problem admits an optimal primal--dual solution pair, which is standard in literature \citep{zhang2022acgd}:
\begin{assumption}\label{ass:pridual_slater}
There exist $x^*\in X$ and
$\lambda^*\in\mathbb R_+^m$ such that
\begin{equation} \label{eq:standing_primal_dual_pair}
    \begin{aligned}
        & \langle
F(x^*)+J(x^*)^\top\lambda^*,
x-x^*
 \rangle
\geq0,
\;\forall x\in X,
 \quad  G(x^*)\leq0, \quad
\langle\lambda^*,G(x^*)\rangle=0.
    \end{aligned}
\end{equation}
In addition, Slater's condition holds: there exists
$x^{\mathrm{sl}}\in ri(X)$ such that
$g_i(x^{\mathrm{sl}})<0$,  $\forall i=1,\ldots,m$.
\end{assumption}

\subsection{Related Literature and Research Objective}
First-order methods for VIs defined over a simple set $X$ are well studied.
For Lipschitz continuous monotone operators $F$, mirror-prox
\citep{nemirovski2004prox} and dual extrapolation \citep{nesterov2007dual} compute an $\varepsilon$-approximate weak solution using $\mathcal O(1/\varepsilon)$ operator evaluations.
This dependence on $\varepsilon$ is optimal in the worst case, matching the $\Omega(1/\varepsilon)$ lower bound established by \citet{ouyang2021lower} for a saddle-point subclass of monotone VIs.

For VIs with functional constraints, however, several research gaps remain in the existing literature.
One main challenge is to design first-order methods that avoid repeated projections onto the functionally constrained set $\widetilde X$.
So far, there are two closely related studies partially closing this gap.
\citet{boob2026firstorder} develop adaptive operator extrapolation methods that, in the smooth setting, achieve $\varepsilon$-accuracy in both the weak VI gap and constraint violation using $\mathcal O(1/\varepsilon)$ evaluations of both the operator and the constraint gradients.
Alternatively, \citet{zhang2025primal} develop a primal gradient-descent method using local quadratic programs.
Without requiring Lipschitz continuity of the operator or information about optimal multipliers, their method achieves these accuracy guarantees in $\mathcal O(1/\varepsilon^2)$ iterations.
This leaves a fundamental question: \emph{can we design first-order methods that require fewer oracle evaluations and provide sharper constraint-violation bounds without relying on local quadratic programs?}

Our work is motivated by results from smooth convex optimization, which
suggest that faster constraint-gradient bounds may be attainable.
For an $\mathcal{O}(1)$ number of constraints, \citet{xu2022fewconstraints} obtain $\widetilde{\mathcal O}(1/\sqrt{\varepsilon})$ first-order evaluation complexity for finding $\varepsilon$-approximate Karush–Kuhn–Tucker (KKT) points, where $\widetilde{\mathcal O}$ suppresses logarithmic factors.
For  the multiple-constraint case, \citet{zhang2022acgd} achieve $\mathcal O(1/\sqrt{\varepsilon})$ evaluation complexity for objective accuracy and feasibility.
These optimization results suggest that
$\mathcal O(1/\sqrt{\varepsilon})$ constraint-gradient
complexity may also be attainable for monotone CVIs.
However, the worst-case $\Omega(1/\varepsilon)$ operator
lower bound discussed above applies even without functional
constraints.
Achieving the desired $\mathcal{O}(1/\sqrt{\varepsilon})$ constraint-gradient bound therefore
requires separating the two oracle costs.
 Gradient sliding \citep{lan2016gradient} provides a way to
achieve such separation by reusing gradient information
across multiple updates involving another component.
For structured VIs whose operator is the sum of a smooth convex gradient and a Lipschitz continuous monotone component, \citet{lan2021mirrorproxsliding} obtain $\mathcal O(1/\sqrt{\varepsilon})$ gradient evaluations and $\mathcal O(1/\varepsilon)$ evaluations of the monotone component.
However,  in CVIs,  the smooth functions define the feasible set rather than an additive component of the operator. Extending this separation to CVIs therefore requires
controlling constraint violation together with the weak VI gap.

Another feature we seek is  parameter-free , meaning that the algorithm can adapt to unknown  Lipschitz constant of the operator and smoothness constants of the constraint functions while retaining the same complexity bounds.
The adaptive golden ratio algorithm \citep{malitsky2020golden}
and past-extragradient methods \citep{ene2022adaptivevi}
compute $\varepsilon$-approximate weak solutions to Lipschitz
continuous monotone VIs on bounded domains using
$\mathcal O(1/\varepsilon)$ operator evaluations,
without prior knowledge of the Lipschitz constant.
However, these methods require proximal updates over the feasible domain, so applying them directly to CVIs retains the difficulty of projecting onto $\widetilde X$.
Adaptive extensions are also available for convex composite optimization \citep{wu2026universal} and convex optimization with functional constraints \citep{zhang2022acgd}.

In summary, our objective is to design 
first-order algorithms that compute an $\varepsilon$-approximate
weak solution to monotone CVIs with convex smooth functional constraints using
$\mathcal O(1/\varepsilon)$ operator evaluations and
$\mathcal O(1/\sqrt{\varepsilon})$ constraint-gradient
evaluations,
without projections onto $\widetilde X$ or prior
knowledge of Lipschitz constant $M$, smoothness constant $L$, or a uniform bound on $\|J(x)\|$ over $X$.

\subsection{Organization and Main Results}
The remainder of this paper is organized as follows.
Section~\ref{sec:pfaos} introduces the parameter-free accelerated
operator sliding method (PFAOS) for the single-constraint case.
PFAOS achieves $\mathcal O(1/\varepsilon)$ operator 
and $\mathcal O(1/\sqrt{\varepsilon})$ constraint-gradient evaluation 
complexities without prior knowledge of the Lipschitz constant
of $F$ or the smoothness constant of $G$ or projections onto $\widetilde X$.
These bounds retain the optimal operator-complexity order for unconstrained VIs while improving the existing constraint-gradient complexity of CVIs from $\mathcal O(1/\varepsilon)$  to $\mathcal O(1/\sqrt{\varepsilon})$.
Section~\ref{sec:mc_pfaos} extends PFAOS to the multi-constrained
parameter-free accelerated operator sliding method (MC-PFAOS)
for solving the multi-constraint case.
MC-PFAOS attains $\mathcal O(1/\varepsilon)$ operator and
$\mathcal O(1/\sqrt{\varepsilon})$ Jacobian evaluation complexities
without projections onto $\widetilde X$ or prior knowledge of
the Lipschitz or smoothness constants or a uniform bound on
$\|J(x)\|$ over $X$.
To the best of our knowledge, these results substantially improve the complexity bounds established for CVI algorithms, including those assuming known smoothness parameters \citep{boob2026firstorder, zhang2025primal}.

Next, we provide two complementary results that relax the requirements needed to implement PFAOS and MC-PFAOS. 
Section~\ref{sec:special_first_iteration} modifies the first-iteration search for both methods to allow an arbitrarily small  initial estimate while preserving the same optimal complexities.
Section~\ref{sec:stopping_criteria} develops an automatic
stopping rule for both methods that requires no knowledge
of the diameter of the set $X$ while preserving these complexities.
Section~\ref{sec:numerical_experiments} reports numerical
experiments on large-scale games with coupling constraints,
and Section~\ref{sec:Conclusion} concludes the paper.

\section{Parameter-free Accelerated Operator Sliding (PFAOS)}\label{sec:pfaos}
Algorithm \ref{alg:pfaos} presents the PFAOS method for solving the monotone CVI  in \eqref{eq:problem_of_interest}  with a single constraint ($m=1$). For notational simplicity, we write $g:=g_1$, throughout this section.

\subsection{Overview of PFAOS}

\begin{algorithm}[!htb]
\caption{Parameter-Free Accelerated Operator Sliding (PFAOS)}
\label{alg:pfaos}
\footnotesize
\begin{algorithmic}[1]
\State{\textbf{Initialization:} $\xub_0 = x_0\in X$, $L_0\in(0,L]$, \(\Lambda>0\), $M_0=\eta_0^{S_0}=2 \Lambda L_0$}
\For{$k=1,\dots,N$}
    \State {Set $s=1$, $M_k^0 = M_{k-1}$}
    \State{Set \label{line:pfaos_s_loop} $\eta_k^s =  2^{s-2}\eta_{k-1}^{S_{k-1}}$
    and compute $\xlb_k^s = (1-\gamma_k^s)\xub_{k-1}+\gamma_k^sx_{k-1}$, where  \label{alg1_line:4}
    \begin{equation}\label{eq:recursion_def}
        \gamma_1^s=1,  \text{ and }\gamma_k^s>0 \text{ satisfying }
        \eta_k^s\gamma_k^s = \eta_{k-1}^{S_{k-1}}\gamma_{k-1}^{S_{k-1}}(1-\gamma_k^s), \forall k\ge2
    \end{equation}
    }
    \State{Set $t=1$,
$M_k^{0,s}=\max\{\eta_k^s,M_k^{s-1}\}$,
$T_k^{0,s}=M_k^{0,s}/\eta_k^s$,
$x_k^{0,s}=x_{k-1}$ and
$x_k^{-1,s}=x_{k-1}^{T_{k-1}^{S_{k-1}}-1,S_{k-1}}$. 
Here, 
$x_{0}^{T_{0}^{S_{0}}-1,S_{0}}=x_0$.\label{alg1_line:5}}
    \State Compute  $x_k^{t,s}$ by solving
\label{line:pfaos_t_loop}
\begin{equation}\label{eq:pfaos_subproblem}
    \begin{aligned}
    \min_{x\in X}\quad
    &  \langle
    F(x_k^{t-1,s})
    +\theta_k^{t,s}
       [F(x_k^{t-1,s})-F(x_k^{t-2,s}) ],
    x
     \rangle
    +\tfrac{\eta_k^s}{2}\|x-x_{k-1}\|^2
    +\tfrac{p_k^{t,s}}{2}\|x-x_k^{t-1,s}\|^2 \\
    \text{s.t.}\quad
    & g(\xlb_k^s)
    + \langle
    \nabla g(\xlb_k^s),x-\xlb_k^s
     \rangle
    \le 0 .
    \end{aligned}
    \end{equation}
    Here $\theta_1^{1,s}=1$, 
    $\theta_k^{1,s}=M_k^{1,s}/M_{k-1}$,
    $\forall k\ge2$, $\theta_k^{t,s}=M_k^{t,s}/M_k^{t-1,s}$,
    $\forall t\ge2$, and $p_k^{t,s} = 2M_k^{t,s}+\eta_k^sM_k^{t,s}\sum\nolimits_{i=1}^{t-1}(1/M_k^{i,s})$,  where $M_k^{t,s} = 2^{q_k^{t,s}}M_k^{t-1,s}$ and $q_k^{t,s}\ge 0$ is the smallest integer such that $\|F(x_k^{t,s})-F(x_k^{t-1,s})\|\le M_k^{t,s}\|x_k^{t,s}-x_k^{t-1,s}\|$.
    \State{Set $T_k^{t,s}\gets t-1+2^{q_k^{t,s}}(T_k^{t-1,s}-t+1)$}\label{line:pfaos_t_termination}
    \State{\textbf{if } $t<T_k^{t,s}$, \textbf{then } set $t\gets t+1$ and go to line \ref{line:pfaos_t_loop}; otherwise let $T_k^s = T_k^{t,s}$, $M_k^s = M_k^{T_k^s,s}$, \label{alg1_line:8}
    \begin{equation*}
            x_k^s = x_k^{T_k^s,s}
            ,\quad 
            \xhat_k^s = (\sum\nolimits_{i=1}^{T_k^s}1/M_k^{i,s})^{-1}\sum\nolimits_{i=1}^{T_k^s}x_k^{i,s}/M_k^{i,s}
            ,\quad
            \xub_k^s = (1-\gamma_k^s)\xub_{k-1}+\gamma_k^s\xhat_k^s
    \end{equation*}
    }
    \State{\textbf{if } $g(\xub_k^s)\le g(\xlb_k^s)+\langle\nabla g(\xlb_k^s), \xub_k^s-\xlb_k^s\rangle+\tfrac{\eta_k^s}{2\Lambda\gamma_k^s}\|\xub_k^s-\xlb_k^s\|^2$ \text{does not hold, } \textbf{then } update $s\gets s+1$ and go to line \ref{line:pfaos_s_loop}; \textbf{else } set $S_k=s$, $M_k=M_k^{S_k}$
    $x_k=x_k^{S_k}$, $\xhat_k = \xhat_k^{S_k}$ and $\xub_k = \xub_k^{S_k}$.}\label{line:pfaos_g_test}
\EndFor

\State{\Return $\xub_N$}
\end{algorithmic}
\end{algorithm}

PFAOS has three iteration levels.
The outer loop, indexed by $k$, generates the accelerated sequence $\{\xub_k\}$.
At each outer iteration, the middle loop, indexed by $s$,   performs a line search for the constraint-curvature parameter $\eta_k^s$.
Each middle trial $(k,s)$ runs an inner operator-sliding loop, indexed by $t$, and updates the local Lipschitz estimate $M_k^{t,s}$ of $F$.
We use the superscript $s$ to index middle-loop trials and $(t,s)$ to index inner-loop iterations within each trial, with $S_k$ denoting the accepted middle-trial index and $T_k^s$ the final inner-loop index of trial $(k,s)$.

 \subsubsection{Constraint Linearization and Operator Sliding}
Operator sliding is implemented  within each middle-trial  $(k,s)$.
Given $\eta_k^s$, the algorithm determines $\gamma_k^s$ using \eqref{eq:recursion_def} and forms $\xlb_k^s$ in line~\ref{alg1_line:4}.
In  subproblem~\eqref{eq:pfaos_subproblem}, we approximate the functional constraint $g$ by its affine linearization with $\nabla g(\xlb_k^s)$.
This constraint remains fixed throughout the corresponding inner loop.

At each inner iteration $t$, the algorithm solves subproblem~\eqref{eq:pfaos_subproblem}, whose objective
combines an extrapolated evaluation of $F$ with two quadratic regularization terms.
Each inner line-search trial requires one new evaluation of $F$,
while the fixed affine constraint is reused without reevaluating
$\nabla g$.
After $T_k^s$ inner iterations, line~\ref{alg1_line:8} retains the last iterate
$x_k^s=x_k^{T_k^s,s}$, computes the weighted average $\xhat_k^s$
with weights proportional to $1/M_k^{t,s}$, and forms the candidate
outer iterate $\xub_k^s$.

\subsubsection{Parameter Adaptation}
PFAOS employs two line searches to adapt the local constraint-curvature parameter $\eta_k^s$ and
the local Lipschitz estimate $M_k^{t,s}$ of $F$, respectively.
For initialization,  $L_0\in(0, L]$ can be obtained from any distinct $u,v\in X$ with positive linearization error by setting
\(
L_0=
2[g(v)-g(u)
-\langle\nabla g(u),v-u\rangle]
/
\|v-u\|^2
.
\)
The smoothness of $g$ guarantees $L_0\le L$, so this construction does not require knowing $L$ a priori.

At the start of middle trial $(k,s)$, the algorithm sets $\eta_k^s =2^{s-2}\eta_{k-1}^{S_{k-1}}.$
We treat the first and subsequent trials separately: the first trial uses half of the previously accepted
estimate, and each subsequent trial doubles the current estimate.
Within each middle trial, the initial estimate $M_k^{0,s}$ of the Lipschitz constant of $F$ is
updated in line~\ref{alg1_line:5}, where $M_k^0=M_{k-1}$.
At inner iteration $t$, the algorithm first tests the preceding estimate $M_k^{t-1,s}$.
For each trial value of $M_k^{t,s}$, it computes the parameters in line \ref{line:pfaos_t_loop},
solves subproblem ~\eqref{eq:pfaos_subproblem}, and tests the condition
\(
\|F(x_k^{t,s})-F(x_k^{t-1,s})\|
\le
M_k^{t,s}\|x_k^{t,s}-x_k^{t-1,s}\|
\).
If this condition fails, the estimate is doubled, the parameters are updated, and the subproblem
is solved again.
Thus, $q_k^{t,s}$ records the smallest number of doublings required to satisfy this condition, with
$M_k^{t,s}=2^{q_k^{t,s}}M_k^{t-1,s}$.

The algorithm  initializes the inner-loop horizon   as $T_k^{0,s}=M_k^{0,s}/\eta_k^s$.
After accepting $M_k^{t,s}$, the algorithm updates $T_k^{t,s}$ in line~\ref{line:pfaos_t_termination}.
This update leaves the first $t-1$ iterations unchanged and scales the remaining horizon by
$M_k^{t,s}/M_k^{t-1,s}$.
In particular, the update preserves
\begin{equation*}
\sum\nolimits_{i=1}^{t}\tfrac{1}{M_k^{i,s}} +\tfrac{T_k^{t,s}-t}{M_k^{t,s}}
=\tfrac{1}{\eta_k^s}.
\end{equation*}
Consequently, a larger operator estimate requires additional inner updates to compensate for the
smaller weights $1/M_k^{t,s}$.
The inner loop terminates when $t=T_k^{t,s}$, and its final horizon is recorded as $T_k^s$.
At termination, the accumulated weights satisfy 
$\sum\nolimits_{t=1}^{T_k^s}1/M_k^{t,s} =1/\eta_k^s.$

After completing the inner loop, PFAOS
checks whether the candidate $\xub_k^s$ satisfies
the quadratic upper-bound condition
\(
g(\xub_k^s)\le g(\xlb_k^s)+\langle\nabla g(\xlb_k^s),\xub_k^s-\xlb_k^s\rangle +\tfrac{\eta_k^s}{2\Lambda\gamma_k^s}\|\xub_k^s-\xlb_k^s\|^2.
\)
If this condition fails, the algorithm goes back to line \ref{line:pfaos_s_loop}, increases $\eta_k^s$, updates the constraint linearization,
and restarts the inner loop from the preceding accepted outer iterate.
Otherwise, the trial is accepted, and its candidate iterates and parameter estimates are retained for outer iteration $k$.

Finally, PFAOS verifies that the adaptive inner loop terminates in finitely many iterations.
The curvature and operator estimates remain dyadic multiples of their common initial value
$M_0=\eta_0^{S_0}$.
Since $M_k^{0,s}\ge\eta_k^s$, the initial horizon $T_k^{0,s}=M_k^{0,s}/\eta_k^s$ is a positive integer.
The horizon update preserves integrality and ensures $T_k^{t,s}\ge t$ throughout the inner loop.
Moreover, the local Lipschitz condition holds whenever $M_k^{t,s}\ge M$, and the minimality of $q_k^{t,s}$ implies $M_k^{t,s}<2M$.
Thus, the operator estimate and inner-loop horizon can increase only finitely many times.

The adaptive inner-loop horizon allows multiple operator updates
per constraint linearization.
Combined with acceleration, this separation is a key to  achieving
$\mathcal O(1/\sqrt{\varepsilon})$ constraint-gradient evaluations
and $\mathcal O(1/\varepsilon)$ operator evaluations, as established
in the following analysis.

\subsection{Convergence Analysis of PFAOS}\label{sec:pfaos_convergence}
In this subsection, we establish the convergence properties of
Algorithm~\ref{alg:pfaos}. 
The following quantities are introduced for the analyses.
For each middle trial $(k,s)$, define the total averaging weight,
the affine constraint model, and the operator differences by
\begin{equation*}
\begin{aligned}
\mathcal B_k^s&:=\sum\nolimits_{t=1}^{T_k^s}\tfrac{1}{M_k^{t,s}}, \quad
\ell_k^s(z) :=g(\xlb_k^s) +\langle\nabla g(\xlb_k^s),z-\xlb_k^s\rangle, \quad
\Delta_k^{t,s}:=F(x_k^{t,s})-F(x_k^{t-1,s}), 0\leq t\leq T_k^s.
\end{aligned}
\end{equation*}
In particular,
$\Delta_k^{0,s}=F(x_{k-1})-F(x_k^{-1,s})$.

Let $\lambda_k^{t,s}\geq0$ be a  KKT  multiplier associated with
$x_k^{t,s}$ in \eqref{eq:pfaos_subproblem}. Using the same averaging weights as for the primal iterates, define
\begin{equation*}
\hat\lambda_k^s:=\tfrac{1}{\mathcal B_k^s} \sum\nolimits_{t=1}^{T_k^s}\tfrac{\lambda_k^{t,s}}{M_k^{t,s}},
 \quad \hat\lambda_k:=\hat\lambda_k^{S_k},  \quad 
\overline\lambda_0=0, \quad \overline\lambda_k:=(1-\gamma_k^{S_k})\overline\lambda_{k-1} +\gamma_k^{S_k}\hat\lambda_k.
\end{equation*}

Denoting $ \Gamma_k^s:=\eta_k^s\gamma_k^s$, by
\eqref{eq:recursion_def}, we have
\(
    \Gamma_k^s
    =
    \Gamma_{k-1}^{S_{k-1}}
     (1-\gamma_k^s )
\), for $k\geq2$, $s\geq1$. 
Given an initial point $x_0\in X$, let \(D_X:=\max_{x\in X}\|x-x_0\| \). Finally, for any $x\in X$ and $\lambda\ge0$, define the 
primal--dual gap function
\begin{equation*}
    \mathcal G_k(x,\lambda)
    :=\langle F(x),\xub_k-x\rangle
      +\lambda g(\xub_k)-\overline{\lambda}_k g(x).
\end{equation*}

Lemma~\ref{thm:three_point} provides a constrained three-point
inequality used repeatedly in the subsequent analysis.
Its proof is a straightforward
adaptation of the argument in \cite[Lemma~3.5]{lan2020first},
using the Lagrangian and the KKT conditions.
\begin{lemma}
\label{thm:three_point}
Let $\mathcal Z\subseteq\mathbb R^n$ be a nonempty closed convex set, let
$q:\mathcal Z\to\mathbb R$ be convex, and let
$c=(c_1,\ldots,c_m):\mathcal Z\to\mathbb R^m$ have convex components.
Given $z,z'\in\mathcal Z$ and $\mu_1,\mu_2\geq0$, suppose that
\(
    v\in\arg\min_{u\in\mathcal Z}
     \{
        q(u)
        +\tfrac{\mu_1}{2}\|u-z\|^2
        +\tfrac{\mu_2}{2}\|u-z'\|^2
        : c(u)\leq0
     \},
\)
and let $\lambda^+\in\mathbb R_+^m$ be a KKT multiplier
associated with $v$. Then, for every $u\in\mathcal Z$ and $\lambda\in\mathbb R_+^m$,
\begin{equation*}\label{eq:constrained_three_point}
    \begin{aligned}
        &q(v)-q(u)
 +\langle\lambda,c(v)\rangle
 -\langle\lambda^+,c(u)\rangle
 +\tfrac{\mu_1+\mu_2}{2}\|v-u\|^2
\nonumber
\leq
 \tfrac{\mu_1}{2}
  (\|u-z\|^2-\|v-z\|^2 )
 +\tfrac{\mu_2}{2}
  (\|u-z'\|^2-\|v-z'\|^2 )
    \end{aligned}
\end{equation*}
\end{lemma}

A key feature of PFAOS is that the inner-loop horizon is not prescribed in advance. Whenever the operator line search increases $M_k^{t,s}$, the
remaining horizon is adjusted accordingly at line \ref{line:pfaos_t_termination}. 
Lemma \ref{thm:t_loop_termination_invariance} shows that this adaptive rule preserves an exact balancing invariant, which is a key property leading the subsequent convergence analyses.

\begin{lemma}\label{thm:t_loop_termination_invariance}
For every middle trial $(k,s)$ executed by Algorithm~\ref{alg:pfaos}, $\eta_k^s\mathcal B_k^s=1$.
\end{lemma}
\begin{proof}{Proof: }

The horizon update and
$M_k^{t,s}=2^{q_k^{t,s}}M_k^{t-1,s}$ give
\(
\tfrac{1}{M_k^{t,s}}
 +\tfrac{T_k^{t,s}-t}{M_k^{t,s}}
 =\tfrac{T_k^{t-1,s}-t+1}{M_k^{t-1,s}}
\). 
Starting from
$T_k^{0,s}/M_k^{0,s}=1/\eta_k^s$ and summing over the
executed inner iterations yields
\(
\sum\nolimits_{i=1}^{t}\tfrac{1}{M_k^{i,s}}
 +\tfrac{T_k^{t,s}-t}{M_k^{t,s}}
 =\tfrac{1}{\eta_k^s}.
\)
At termination, $t=T_k^{t,s}=T_k^s$, so the remaining-horizon
term vanishes. \hfill \Halmos
\end{proof}

Using the invariant in Lemma~\ref{thm:t_loop_termination_invariance},
Proposition~\ref{thm:pfaos_inner_analysis} establishes an
inner-loop bound for each middle trial $(k,s)$, which will
be used in the subsequent outer-loop analysis.

\begin{proposition}\label{thm:pfaos_inner_analysis}
Fix any $k\geq1$ and 
$1\leq s\leq S_k$. For any point $x\in X$ and  scalar $\lambda\geq0$, the iterates generated by Algorithm \ref{alg:pfaos} guarantee 
\(
    \mathcal B_k^s
     \langle F(x),\xhat_k^s-x \rangle
    \leq
    U_k^s+V_k^s+W_k^s,
\)
where
\begin{equation*}
    \begin{aligned}
        U_k^s:=\;&
 \tfrac{1}{M_k^s}
  \langle
    \Delta_k^{T_k^s,s},
    x_k^s-x
  \rangle
-\tfrac{\theta_k^{1,s}}{M_k^{1,s}}
  \langle
   \Delta_k^{0,s},
    x_{k-1}-x
  \rangle
+\tfrac{p_k^{1,s}}{2M_k^{1,s}}
    \|x-x_{k-1}\|^2
 -\tfrac{\eta_k^s+p_k^{T_k^s,s}}{2M_k^s}
    \|x-x_k^s\|^2,
\\
V_k^s:=\;&
 \tfrac{\theta_k^{1,s}}{M_k^{1,s}} 
  \|\Delta_k^{0,s}\|
  \|x_k^{1,s}-x_{k-1} \|
-\tfrac{p_k^{1,s}}{4M_k^{1,s}}
  \|x_k^{1,s}-x_{k-1} \|^2
-\tfrac{p_k^{T_k^s,s}}{4M_k^s}
  \|x_k^s-x_k^{T_k^s-1,s} \|^2,
\\ 
W_k^s:=\;&
 \mathcal B_k^s\hat{\lambda}_k^s\ell_k^s(x)
-\mathcal B_k^s\lambda\ell_k^s(\xhat_k^s)
+\tfrac{1}{2}
    \|x-x_{k-1}\|^2
 -\tfrac{1}{2}
    \|\xhat_k^s-x_{k-1}\|^2.
    \end{aligned}
\end{equation*}
\end{proposition}
\begin{proof}{Proof:}
Fix any $t\ge1$. Applying Lemma \ref{thm:three_point} to the inner
subproblem \eqref{eq:pfaos_subproblem} yields that
\begin{equation*}
\begin{aligned}
& \langle
    F(x_k^{t-1,s})
    +\theta_k^{t,s}\Delta_k^{t-1,s}
    +\lambda_k^{t,s}\nabla g(\xlb_k^s),
    x_k^{t,s}-x
  \rangle
 + \tfrac{\eta_k^s+p_k^{t,s}}{2}
    \|x-x_k^{t,s}\|^2\\
\leq &
  \tfrac{\eta_k^s}{2}
  (
    \|x-x_{k-1}\|^2
    -\|x_k^{t,s}-x_{k-1}\|^2
  )
 + \tfrac{p_k^{t,s}}{2}
  (
    \|x-x_k^{t-1,s}\|^2
    -\|x_k^{t,s}-x_k^{t-1,s}\|^2
  ).
  \end{aligned}
\end{equation*}
The primal feasibility and complementary slackness of this inner subproblem also guarantee that $-\lambda\ell_k^s(x_k^{t,s})\ge0$ and $\lambda_k^{t,s}\ell_k^s(x_k^{t,s})=0$.
Combining the above relations with the inequality \(
     \langle F(x)-F(x_k^{t,s}),x_k^{t,s}-x \rangle
    \leq
     0
\) coming from the monotonicity of $F$, applying the
Cauchy--Schwarz
inequality, and rearranging terms, we obtain
\(
     \langle F(x),x_k^{t,s}-x \rangle
    \leq
    U_k^{t,s}+V_k^{t,s}+W_k^{t,s},
\)
where
\begin{equation*}
    \begin{aligned}
        U_k^{t,s}:=\;&
     \langle\Delta_k^{t,s},x_k^{t,s}-x \rangle
    -\theta_k^{t,s}
      \langle
        \Delta_k^{t-1,s},x_k^{t-1,s}-x
      \rangle
+ \tfrac{p_k^{t,s}}{2}\|x-x_k^{t-1,s}\|^2
 - \tfrac{\eta_k^s+p_k^{t,s}}{2}
     \|x-x_k^{t,s}\|^2,
\\
V_k^{t,s}:=\;&
    \theta_k^{t,s}
    \|\Delta_k^{t-1,s}\|
    \|x_k^{t,s}-x_k^{t-1,s}\|
- \tfrac{p_k^{t,s}}{2}
    \|x_k^{t,s}-x_k^{t-1,s}\|^2,
\\
W_k^{t,s}:=\;&
    -\lambda\ell_k^s(x_k^{t,s})
    +\lambda_k^{t,s}\ell_k^s(x)
+ \tfrac{\eta_k^s}{2}
  [
    \|x-x_{k-1}\|^2
    -\|x_k^{t,s}-x_{k-1}\|^2
  ].
    \end{aligned}
\end{equation*}
Multiplying both sides of the above relation by
$1/M_k^{t,s}$ and summing over $t=1,\ldots,T_k^s$, while considering the definition of $\xhat_k^s$ yield
\begin{equation*}
    \mathcal B_k^s
 \langle F(x),\xhat_k^s-x \rangle
\leq
\sum\nolimits_{t=1}^{T_k^s} \tfrac{U_k^{t,s}}{M_k^{t,s}}
+\sum\nolimits_{t=1}^{T_k^s} \tfrac{V_k^{t,s}}{M_k^{t,s}}   
+\sum\nolimits_{t=1}^{T_k^s} \tfrac{W_k^{t,s}}{M_k^{t,s}}.
\end{equation*}

We bound the right-hand side in three parts: 
(i)  The settings of $\theta_k^{t,s}$ and $p_k^{t,s}$ in Algorithm \ref{alg:pfaos} give  $  \theta_k^{t,s}/M_k^{t,s}
    =
     1/M_k^{t-1,s}, 
    p_k^{t,s}/M_k^{t,s}
    =
     (\eta_k^s+p_k^{t-1,s})/M_k^{t-1,s}$ for $ 2\leq t\leq T_k^s$. 
Then, canceling the recursive terms while noticing the definitions of $M_k^s$, $x_k^s$ and $x_k^{0,s}$ yields $\sum\nolimits_{t=1}^{T_k^s} U_k^{t,s}/M_k^{t,s}
=U_k^s.$
(ii) The inner line search gives
$\|\Delta_k^{t-1,s}\|\leq M_k^{t-1,s}\|x_k^{t-1,s}-x_k^{t-2,s}\|$ for $t\geq2$.
  $p_k^{t,s}/M_k^{t,s}\geq2$ yields
\begin{equation*}
\begin{aligned}
\tfrac{\theta_k^{t,s}}{M_k^{t,s}}
 \|\Delta_k^{t-1,s}\|\|x_k^{t,s}-x_k^{t-1,s}\| 
 \leq
 \tfrac{p_k^{t-1,s}}{4M_k^{t-1,s}}
       \|x_k^{t-1,s}-x_k^{t-2,s}\|^2
 +\tfrac{p_k^{t,s}}{4M_k^{t,s}}
       \|x_k^{t,s}-x_k^{t-1,s}\|^2.
\end{aligned}
\end{equation*}
Thus, the interior cross terms are absorbed by the negative squared-step terms, leaving the boundary expression $\sum_{t=1}^{T_k^s}V_k^{t,s}/M_k^{t,s}\le V_k^s$.
This also covers $T_k^s=1$, when there are no interior cross terms.
(iii) The affinity of $\ell_k^s$ and the definitions of the weighted averages give the following:
\begin{equation*}
\sum\nolimits_{t=1}^{T_k^s}\tfrac{W_k^{t,s}}{M_k^{t,s}}=
 \mathcal B_k^s\hat\lambda_k^s\ell_k^s(x)  -\mathcal B_k^s\lambda\ell_k^s(\xhat_k^s)
 +\tfrac{\eta_k^s\mathcal B_k^s}{2}\|x-x_{k-1}\|^2 -\tfrac{\eta_k^s}{2}\sum\nolimits_{t=1}^{T_k^s}   \tfrac{\|x_k^{t,s}-x_{k-1}\|^2}{M_k^{t,s}}.
\end{equation*} 
By the convexity of the squared norm, we have $\sum\nolimits_{t=1}^{T_k^s}
 \tfrac{\|x_k^{t,s}-x_{k-1}\|^2}{M_k^{t,s}}
 \geq\mathcal B_k^s\|\xhat_k^s-x_{k-1}\|^2.$
Using the Lemma \ref{thm:t_loop_termination_invariance}, we obtain the bound $\sum\nolimits_{t=1}^{T_k^s}\tfrac{W_k^{t,s}}{M_k^{t,s}}\le W_k^s$. Combining the three bounds proves the result. \hfill
 \Halmos
\end{proof}

Combining Proposition~\ref{thm:pfaos_inner_analysis} at $s=S_k$
with the middle-loop acceptance condition at line \ref{line:pfaos_g_test} of Algorithm \ref{alg:pfaos} yields the following
one-step bound on the outer primal--dual gap.

\begin{proposition}\label{thm:pfaos_middle_gap}
For every $k\geq1$, $x\in X$ and $\lambda\geq0$, the iterates generated by
Algorithm~\ref{alg:pfaos} satisfy
\begin{equation}\label{eq:pfaos_outer_gap}
\mathcal G_k(x,\lambda)
- (1-\gamma_k^{S_k} )\mathcal G_{k-1}(x,\lambda)
\leq H_k+I_k+R_k(\lambda),
\end{equation}
where $R_k(\lambda):=
  \tfrac{\gamma_k^{S_k}\eta_k^{S_k}}{2}
  ( \tfrac{\lambda}{\Lambda}-1 )
 \|\xhat_k-x_{k-1}\|^2$, and 
\begin{equation*} 
    \begin{aligned}
        H_k:=&
  \tfrac{\gamma_k^{S_k}}
      {\mathcal B_k^{S_k}M_k^{S_k}}
  \langle
   \Delta_k^{T_k^{S_k},S_k},
   x_k-x
  \rangle
-
  \tfrac{\gamma_k^{S_k}\theta_k^{1,S_k}}
      {\mathcal B_k^{S_k}M_k^{1,S_k}}
  \langle
   \Delta_k^{0,S_k},
   x_{k-1}-x
  \rangle
\\
&+
  (
    \tfrac{\gamma_k^{S_k}}{2\mathcal{B}_k^{S_k}}
   +
    \tfrac{\gamma_k^{S_k}p_k^{1,S_k}}
        {2\mathcal B_k^{S_k}M_k^{1,S_k}}
  )
 \|x-x_{k-1}\|^2
-
  \tfrac{
   \gamma_k^{S_k}
    (\eta_k^{S_k}
         +p_k^{T_k^{S_k},S_k} )}
 {2\mathcal B_k^{S_k}M_k^{S_k}}
 \|x-x_k\|^2,
\\
I_k:=&
  \tfrac{\gamma_k^{S_k}\theta_k^{1,S_k}}
      {\mathcal B_k^{S_k}M_k^{1,S_k}}
 \|\Delta_k^{0,S_k}\|
 \|x_k^{1,S_k}-x_{k-1}\|
-
  \tfrac{\gamma_k^{S_k}p_k^{1,S_k}}
      {4\mathcal B_k^{S_k}M_k^{1,S_k}}
 \|x_k^{1,S_k}-x_{k-1}\|^2
-
  \tfrac{\gamma_k^{S_k}p_k^{T_k^{S_k},S_k}}
      {4\mathcal B_k^{S_k}M_k^{S_k}}
 \|x_k
      -x_k^{T_k^{S_k}-1,S_k}\|^2.
    \end{aligned}
\end{equation*}
\end{proposition}

\begin{proof}{Proof:}
From the definitions of $\xub_k$ and $\overline{\lambda}_k$, it follows that
\begin{equation*}
    \begin{aligned}
        &\mathcal G_k(x,\lambda)
- (1-\gamma_k^{S_k} )\mathcal G_{k-1}(x,\lambda)
=\;
\gamma_k^{S_k}
  \langle F(x),\xhat_k^{S_k}-x \rangle
+\lambda g(\xub_k)
- (1-\gamma_k^{S_k} )\lambda g(\xub_{k-1})
-\gamma_k^{S_k}\hat{\lambda}_k^{S_k}g(x).
    \end{aligned}
\end{equation*}
Here, we make two observations.
First,  Proposition~\ref{thm:pfaos_inner_analysis} and    $\gamma_k^{S_k}/\mathcal B_k^{S_k}=\Gamma_k^{S_k}$ give $ \gamma_k^{S_k}  \langle F(x),\xhat_k^{S_k}-x \rangle \leq \Gamma_k^{S_k} (U_k^{S_k}+V_k^{S_k}+W_k^{S_k} ).$
 
Second, the acceptance test described at line \ref{line:pfaos_g_test}  together with the definitions of $\xlb_k^{S_k}$ and $\xub_k$ ensures
\begin{equation*}
g(\xub_k)
\leq
\ell_k^{S_k}(\xub_k)
+
 \tfrac{\eta_k^{S_k}}
     {2\Lambda\gamma_k^{S_k}}
\|\xub_k-\xlb_k^{S_k}\|^2
=
\ell_k^{S_k}(\xub_k)
+
 \tfrac{\eta_k^{S_k}\gamma_k^{S_k}}
     {2\Lambda}
\|x_{k-1}-\xhat_k\|^2,
\end{equation*}
and 
the convexity of $g$ gives that 
$-g(\xub_{k-1})\le-
\ell_k^{S_k}(\xub_{k-1})$ and $-g(x)\le -\ell_k^{S_k}(x)$. Then, it follows 
\begin{equation*}
\begin{aligned}
\lambda g(\xub_{k}) -(1-\gamma_k^{S_k})\lambda g(\xub_{k-1})
\le
\tfrac{\eta_k^{S_k}\gamma_k^{S_k}\lambda}
     {2\Lambda}
\|x_{k-1}-\xhat_k\|^2 
+\lambda \gamma_k^{S_k}\ell_k^{S_k}(\xhat_k).
\end{aligned}
\end{equation*}
Then, the conclusion follows immediately  with   Lemma \ref{thm:t_loop_termination_invariance} and definitions of $x_k$ and $\xhat_k$. \hfill \Halmos
\end{proof}

Using the conclusion of Propositions~\ref{thm:pfaos_middle_gap}, Theorem~\ref{thm:pfaos_global_convergence} establishes a global
primal--dual gap bound.

\begin{theorem}\label{thm:pfaos_global_convergence}
For every $N\geq1$, $x\in X$, and $\lambda\geq0$, the iterates generated by
Algorithm~\ref{alg:pfaos} satisfy
\begin{equation}\label{eq:pfaos_global_gap}
\mathcal G_N(x,\lambda)
\leq
\tfrac{\Gamma_N^{S_N}}{2}
  (
    3\|x-x_0\|^2
    +
      (\tfrac{\lambda}{\Lambda}-1  )
    \sum\nolimits_{k=1}^N
    \|\xhat_k^{S_k}-x_{k-1}\|^2
  ).
\end{equation}
In particular, if  $\Lambda\ge \lambda$,
\(
    \mathcal G_N(x,\lambda)
    \leq
    (3\Gamma_N^{S_N}/2)
    \|x-x_0\|^2.
\)
\end{theorem}
\begin{proof}{Proof:}
Notice that the conclusion of Proposition \ref{thm:pfaos_middle_gap} holds for all $k\ge1$, $x\in X$ and $\lambda\ge0$.
Dividing both sides of \eqref{eq:pfaos_outer_gap} by $\Gamma_k^{S_k}$ for each $k\ge1$, then summing from $k=1$ to $N$ while considering the definition and recursive relation of $\Gamma_k^{S_k}$ yields
\begin{equation}\label{eq:pfaos_outer_gap_sum}
\tfrac{\mathcal G_N(x,\lambda)}{\Gamma_N^{S_N}}
\leq
\sum\nolimits_{k=1}^N \tfrac{H_k}{\Gamma_k^{S_k}}  
+
\sum\nolimits_{k=1}^N \tfrac{I_k}{\Gamma_k^{S_k}} 
+
\sum\nolimits_{k=1}^N \tfrac{R_k(\lambda)}{\Gamma_k^{S_k}}.
\end{equation}
Here, we make three observations. 
First, the initialization of the inner iterates gives $\Delta_{k+1}^{0,S_{k+1}}=\Delta_k^{T_k^{S_k},S_k}$.
Together with $\Delta_1^{0,S_1}=0$, this makes the normalized $H_k$ terms telescope to 
\begin{equation*}
\sum\nolimits_{k=1}^{N}\tfrac{H_k}{\Gamma_k^{S_k}}
 =\tfrac{1}{M_N}\langle\Delta_N^{T_N^{S_N},S_N},x_N-x\rangle
  +\tfrac32\|x-x_0\|^2-\tfrac32\|x-x_N\|^2.
\end{equation*}

Second, for $k\ge2$, the final inner-loop acceptance condition at outer iteration $k-1$, together with the relation 
$\Delta_k^{0,S_k} = \Delta_{k-1}^{T_{k-1}^{S_{k-1}},S_{k-1}}$, gives  $ \|\Delta_k^{0,S_k}\| 
 \leq M_{k-1} \|x_{k-1}-x_{k-1}^{T_{k-1}^{S_{k-1}}-1,S_{k-1}}\|$, which leads 
\begin{equation*}
\begin{aligned}
&\tfrac{\|\Delta_k^{0,S_k}\|}{M_{k-1}}
       \|x_k^{1,S_k}-x_{k-1}\| \leq\tfrac12\|x_k^{1,S_k}-x_{k-1}\|^2
 +\tfrac12\|x_{k-1}
   -x_{k-1}^{T_{k-1}^{S_{k-1}}-1,S_{k-1}}\|^2.
\end{aligned}
\end{equation*}
Absorbing the first term by $I_k/\Gamma_k^{S_k}$
and the second  term by $I_{k-1}/\Gamma_{k-1}^{S_{k-1}}$ 
leaves
\begin{equation*}
\sum\nolimits_{k=1}^{N}\tfrac{I_k}{\Gamma_k^{S_k}}
 \leq-\tfrac{p_N^{T_N^{S_N},S_N}}{4M_N}
       \|x_N-x_N^{T_N^{S_N}-1,S_N}\|^2.
\end{equation*}

Finally, the final inner-loop acceptance condition at the accepted trial
$s=S_N$ gives
\(
\|
\Delta_N^{T_N^{S_N},S_N}
\|
\leq
M_N^{S_N}
\|
x_N-x_N^{T_N^{S_N}-1,S_N}
\|
\).
Moreover, since
\(
p_N^{T_N^{S_N},S_N}/M_N^{S_N}\geq2,
\)
the Cauchy--Schwarz inequality yields
\begin{equation*}
    \tfrac{1}{M_N^{S_N}}
 \langle
\Delta_N^{T_N^{S_N},S_N},
x_N-x
 \rangle
-\tfrac32\|x-x_N\|^2
-(p_N^{T_N^{S_N},S_N}/4M_N^{S_N})
 \|
x_N-x_N^{T_N^{S_N}-1,S_N}
 \|^2
\leq 0.
\end{equation*}
Substituting into \eqref{eq:pfaos_outer_gap_sum} and using the definition of $R_k(\lambda)$ proves \eqref{eq:pfaos_global_gap}.
\hfill \Halmos
\end{proof}

\subsection{Convergence Rate and Evaluation Complexity of PFAOS}\label{sec:pfaos_complexity}
The preceding subsection provides an upper bound of the duality gap involving a scalar $\Gamma_N^{S_N}$. 
In this subsection, we first estimate $\Gamma_k^{S_k}$ and then show the evaluation complexities of PFAOS.
\begin{lemma}\label{thm:pfaos_Gamma_rate}
Every middle-loop line search  in Algorithm~\ref{alg:pfaos}  terminates after finitely many trials, namely, $S_k<\infty$ is well defined and that
\(
    \eta_k^{S_k}/\gamma_k^{S_k}
    \leq 4\Lambda L
\)
and
\(
    \Gamma_k^{S_k}
    \leq
    16\Lambda L/(k+1)^2
\), $\forall k\geq1$.
\end{lemma}

\begin{proof}{Proof:}
Given the $L$-smoothness of $g$, the middle-loop acceptance test in line \ref{line:pfaos_g_test} is
guaranteed whenever
\(
    \eta_k^s/\gamma_k^s\geq\Lambda L
\).
We first verify that    $0<\gamma_k^s\leq1$, for every $(k,s)$ trial, and the middle line search terminates within finitely many trials, namely $S_k<\infty$ is well defined.
For $k=1$, Algorithm~\ref{alg:pfaos} sets $\gamma_1^s=1$, hence 
\(
    \eta_1^s/\gamma_1^s=\eta_1^s.
\)
Since $\eta_1^s$ is doubled after every rejection by line \ref{line:pfaos_g_test}, it eventually satisfies
the preceding sufficient acceptance condition. Thus, we have $S_1<\infty$, and
\(
    \Gamma_1^{S_1}
    =\eta_1^{S_1}\gamma_1^{S_1}>0.
\)
Suppose that $S_{k-1}<\infty$ and $\Gamma_{k-1}^{S_{k-1}}>0$ for some $k\ge2$.  Then the recursion \eqref{eq:recursion_def} gives
\(
    \gamma_k^s
    =
    \Gamma_{k-1}^{S_{k-1}}/
         (\Gamma_{k-1}^{S_{k-1}}+\eta_k^s)
         \in (0,1)
\). 
Indeed, if the middle line search did not terminate, then Algorithm~\ref{alg:pfaos} would generate
\(
    \eta_k^s
    =
    2^{s-2}\eta_{k-1}^{S_{k-1}}
    \to +\infty
    \text{ as }s\to\infty.
\)
Since $\Lambda L<\infty$, there would then exist a finite trial $s$ such
that
\(
  \eta_k^s/\gamma_k^s
    \geq\eta_k^s
    \geq\Lambda L,
\)
at which the middle-loop acceptance test is guaranteed to hold. Thus, $S_k<\infty$ is well defined for all $ k\ge1$.

Next, we prove \(
    \eta_k^{S_k}/\gamma_k^{S_k}
    \leq 4\Lambda L
\).
For $k=1$, we have $\gamma_1^s=1$. If $S_1=1$, then
\(
    \eta_1^{S_1}/\gamma_1^{S_1}
    =\eta_0^{S_0}/2
    =\Lambda L_0
    \le\Lambda L.
\)
Otherwise, the trial immediately preceding the accepted trial failed and
therefore satisfies
\(
    \eta_1^{S_1-1}/\gamma_1^{S_1-1}
    =\eta_1^{S_1-1}
    <\Lambda L.
\)
Since the next trial doubles the previous $\eta_1^s$, the accepted value satisfies
\(
    \eta_1^{S_1}/\gamma_1^{S_1}
    =2\eta_1^{S_1-1}
    <2\Lambda L.
\)
For $k\geq2$ and $s\ge 1$, the recursion \eqref{eq:recursion_def} gives
\(
    \eta_k^s/\gamma_k^s
    =
    \eta_k^s
    +(\eta_k^s)^2/\Gamma_{k-1}^{S_{k-1}},
\) and hence
for the first middle-loop trial at iteration $k$, 
we have
\begin{equation*}
    \tfrac{\eta_k^1}{\gamma_k^1}
    =
    \tfrac{\eta_{k-1}^{S_{k-1}}}{2}
    +
    \tfrac{(\eta_{k-1}^{S_{k-1}})^2}
         {4\eta_{k-1}^{S_{k-1}}
    \gamma_{k-1}^{S_{k-1}}} 
    =
     (
        \tfrac14+\tfrac{\gamma_{k-1}^{S_{k-1}}}{2}
     )
    \tfrac{\eta_{k-1}^{S_{k-1}}}
         {\gamma_{k-1}^{S_{k-1}}}
    \leq
    \tfrac{3\eta_{k-1}^{S_{k-1}}}
         {4\gamma_{k-1}^{S_{k-1}}}.
\end{equation*}
Therefore, if $S_k=1$, then
the induction hypothesis yields
\(
    \eta_k^{S_k}/\gamma_k^{S_k}
    \leq 3\Lambda L
    \leq 4\Lambda L
\). 
Otherwise, $S_k\geq2$ and 
\(
    \eta_k^{S_k-1}/\gamma_k^{S_k-1}
    <\Lambda L
\). 
Since $\eta_k^s=2\eta_k^{s-1}$, we have
\begin{equation*}
     \tfrac{\eta_k^{S_k}}{\gamma_k^{S_k}}
    =
    2\eta_k^{S_k-1}
    + \tfrac{4(\eta_k^{S_k-1})^2}
          {\Gamma_{k-1}^{S_{k-1}}}
    =
    4(
        \eta_k^{S_k-1}
        + \tfrac{(\eta_k^{S_k-1})^2}
              {\Gamma_{k-1}^{S_{k-1}}}
      )
      -2\eta_k^{S_k-1}
    =
    4 \tfrac{\eta_k^{S_k-1}}{\gamma_k^{S_k-1}}
      -2\eta_k^{S_k-1}
    <
    4 \tfrac{\eta_k^{S_k-1}}{\gamma_k^{S_k-1}}
    <4\Lambda L,
\end{equation*}
which proves by induction that
\(
    \eta_k^{S_k}/\gamma_k^{S_k}
    \leq4\Lambda L
\), for $k\geq1$.

It remains to bound $\Gamma_k^{S_k}$. For $k=1$, we have
\(
    \Gamma_1^{S_1}
    =
    \eta_1^{S_1}\gamma_1^{S_1}
    =
    \eta_1^{S_1}/\gamma_1^{S_1}
    \leq4\Lambda L
\). For $k\geq2$, from the recursion \eqref{eq:recursion_def} and the fact 
$1-\sqrt{1-\gamma}\geq\gamma/2$ for $\gamma\in[0,1]$, it follows 
\begin{equation*}
    \tfrac{1}{\sqrt{\Gamma_k^{S_k}}}
    -\tfrac{1}{\sqrt{\Gamma_{k-1}^{S_{k-1}}}}
    =
    \tfrac{1-\sqrt{1-\gamma_k^{S_k}}}
         {\sqrt{\Gamma_k^{S_k}}} 
    \geq
    \tfrac{\gamma_k^{S_k}}
         {2\sqrt{\Gamma_k^{S_k}}} 
    =
    \tfrac{1}
         {2\sqrt{\eta_k^{S_k}/\gamma_k^{S_k}}} 
    \geq
    \tfrac{1}{4\sqrt{\Lambda L}}.
\end{equation*}
Summing gives
\begin{equation*}
     \tfrac{1}{\sqrt{\Gamma_N^{S_N}}} 
    = \tfrac{1}{\sqrt{\Gamma_1^{S_1}}} + \sum\nolimits_{k=2}^N (\tfrac{1}{\sqrt{\Gamma_k^{S_k}}}
    -\tfrac{1}{\sqrt{\Gamma_{k-1}^{S_{k-1}}}})
    \geq
    \tfrac{1}{2\sqrt{\Lambda L}}
    +\tfrac{N-1}{4\sqrt{\Lambda L}}
    =
    \tfrac{N+1}{4\sqrt{\Lambda L}}. \Halmos
\end{equation*} 
\end{proof}

Combining Lemma~\ref{thm:pfaos_Gamma_rate} with Theorem~\ref{thm:pfaos_global_convergence} gives an explicit upper bound for the primal--dual gap. 
We next translate this bound into guarantees for the weak VI gap
and constraint violation in \eqref{eq:problem_of_interest}.

\begin{corollary}\label{thm:pfaos_duality_gap_rate}
Let $(x^*,\lambda^*)$ satisfy \eqref{eq:standing_primal_dual_pair}.
For every $N\geq1$, both
$\sup_{x\in\widetilde X}\langle F(x),\xub_N-x\rangle$ and
$(\Lambda-\lambda^*)[g(\xub_N)]_+$ are bounded above by
\(
\tfrac32\Gamma_N^{S_N}D_X^2
\leq\tfrac{24\Lambda L D_X^2}{(N+1)^2}.
\)
In particular, if $\Lambda\geq\lambda^*+1$, then $\xub_N$ is an
$\varepsilon$-approximate solution of \eqref{eq:problem_of_interest}
whenever
\(
\tfrac32\Gamma_N^{S_N}D_X^2\leq\varepsilon,
\)
which is guaranteed by
\(
N+1\geq\sqrt{\tfrac{24\Lambda L D_X^2}{\varepsilon}}.
\)
\end{corollary}

\begin{proof}{Proof:}
For any $x\in\widetilde X$, since $g(x)\leq0$ and
$\overline\lambda_N\geq0$,  Theorem~\ref{thm:pfaos_global_convergence} gives
\(
\langle F(x),\xub_N-x\rangle
\leq\mathcal G_N(x,0)
\leq\tfrac32\Gamma_N^{S_N}\|x-x_0\|^2.
\)
Taking the supremum over $x\in\widetilde X$ yields the  weak VI gap bound.
Next, by the conditions in  \eqref{eq:standing_primal_dual_pair} and the convexity of $g$, we have 
\(
\langle F(x^*),\xub_N-x^*\rangle
\geq-\lambda^*\langle\nabla g(x^*),\xub_N-x^*\rangle
\geq-\lambda^*(g(\xub_N)-g(x^*))
=-\lambda^*g(\xub_N).
\)
It follows that
\(
(\Lambda-\lambda^*)[g(\xub_N)]_+
\leq[\langle F(x^*),\xub_N-x^*\rangle
+\Lambda g(\xub_N)]_+
\leq[\mathcal G_N(x^*,\Lambda)]_+
\leq\tfrac32\Gamma_N^{S_N}D_X^2.
\)
The explicit rate directly follows from Lemma~\ref{thm:pfaos_Gamma_rate}. \hfill \Halmos
\end{proof}

Corollary~\ref{thm:pfaos_duality_gap_rate} gives a sufficient number
of outer iterations to obtain an $\varepsilon$-approximate solution.
Using this result, Proposition~\ref{thm:pfaos_oracle_complexity} bounds
the numbers of constraint-gradient and operator evaluations.

\begin{proposition}\label{thm:pfaos_oracle_complexity}
Suppose $\Lambda\geq\lambda^*+1$. For any $\varepsilon>0$,
Algorithm~\ref{alg:pfaos} obtains an $\varepsilon$-approximate solution
of \eqref{eq:problem_of_interest} using at most
\(
3+2\sqrt{\tfrac{24\Lambda L D_X^2}{\varepsilon}}
+\lceil\log_2\tfrac{2L}{L_0}\rceil
\)
evaluations of $\nabla g$ and
\(
1+\tfrac{4\overline M}{\eta_0^{S_0}}
+\tfrac{6\overline M D_X^2}{\varepsilon}
+\lceil\log_2\tfrac{\overline M}{M_0}\rceil
\)
evaluations of $F$, where $\overline M:=2\max\{2\Lambda L,M\}$.
These bounds include evaluations at initialization and in
accepted and rejected line-search trials.
\end{proposition}

\begin{proof}{Proof:}
Let $N_\varepsilon\geq1$ be the first outer iteration satisfying
\(
\tfrac32\Gamma_{N_\varepsilon}^{S_{N_\varepsilon}}D_X^2
\leq\varepsilon.
\)
By Corollary~\ref{thm:pfaos_duality_gap_rate},
$\xub_{N_\varepsilon}$ is an $\varepsilon$-approximate solution to \eqref{eq:problem_of_interest} and
\(
N_\varepsilon
\leq1+\sqrt{\tfrac{24\Lambda L D_X^2}{\varepsilon}}
\).
Each middle trial requires one evaluation of $\nabla g$.
The update $\eta_k^{S_k}=2^{S_k-2}\eta_{k-1}^{S_{k-1}}$ gives
\(
\sum\nolimits_{k=1}^{N_\varepsilon}S_k
=2N_\varepsilon
+\log_2\tfrac{\eta_{N_\varepsilon}^{S_{N_\varepsilon}}}{\eta_0^{S_0}}
\leq2N_\varepsilon+\log_2\tfrac{2L}{L_0},
\)
where the inequality follows from
$\eta_{N_\varepsilon}^{S_{N_\varepsilon}}\le\eta_{N_\varepsilon}^{S_{N_\varepsilon}}/\gamma_{N_\varepsilon}^{S_{N_\varepsilon}} \leq4\Lambda L$
by Lemma~\ref{thm:pfaos_Gamma_rate} and $\eta_0^{S_0}=2\Lambda L_0$.
Allowing one initialization evaluation and using the bound on
$N_\varepsilon$ yields the stated constraint-gradient complexity.

For the operator evaluations, note that $M_0\leq2\Lambda L$ and
$\eta_k^s\leq\eta_k^{S_k}\leq4\Lambda L$ for every executed middle
trial. Thus, the initialization
$M_k^{0,s}=\max\{\eta_k^s,M_k^{s-1}\}$ preserves the bound
$\overline M$. Each inner update either leaves the estimate unchanged
or, when backtracking occurs, gives $M_k^{t,s}<2M$ by the
$M$-Lipschitz continuity of $F$ and the minimality of $q_k^{t,s}$.
Induction therefore yields $M_k^{t,s}\leq\overline M$ for every
executed inner iteration. Using this bound and
Lemma~\ref{thm:t_loop_termination_invariance}, we obtain
\(
T_k^s
\leq\sum\nolimits_{t=1}^{T_k^s}\tfrac{\overline M}{M_k^{t,s}}
=\tfrac{\overline M}{\eta_k^s}.
\)
Moreover, $\eta_k^s=2^{s-2}\eta_{k-1}^{S_{k-1}}$ gives
\(
\sum\nolimits_{s=1}^{S_k}\tfrac1{\eta_k^s}
<4/\eta_{k-1}^{S_{k-1}}.
\)
Hence, for $N_\varepsilon\geq2$,
\begin{equation*}
\sum\nolimits_{k=1}^{N_\varepsilon}
\sum\nolimits_{s=1}^{S_k}T_k^s
<4\overline M (
\tfrac1{\eta_0^{S_0}}
+\sum\nolimits_{k=1}^{N_\varepsilon-1}
\tfrac1{\eta_k^{S_k}} )
=4\overline M (
\tfrac1{\eta_0^{S_0}}
+\tfrac1{\Gamma_{N_\varepsilon-1}^{S_{N_\varepsilon-1}}}
 )
<\tfrac{4\overline M}{\eta_0^{S_0}}
+\tfrac{6\overline M D_X^2}{\varepsilon}.
\end{equation*}
The equality follows by telescoping
\(
1/\Gamma_k^{S_k}-1/\Gamma_{k-1}^{S_{k-1}}
=1/\eta_k^{S_k}
\)
for $k\geq2$, with $\Gamma_1^{S_1}=\eta_1^{S_1}$.
The last inequality follows from
\(
\tfrac32\Gamma_{N_\varepsilon-1}^{S_{N_\varepsilon-1}}D_X^2
>\varepsilon
\),
by the minimality of $N_\varepsilon$.
For $N_\varepsilon=1$, the geometric-series bound gives
$\sum_{s=1}^{S_1}T_1^s<4\overline M/\eta_0^{S_0}$.
Finally, each rejected inner trial doubles the nondecreasing operator
estimate, so
\(
\sum\nolimits_{k=1}^{N_\varepsilon}
\sum\nolimits_{s=1}^{S_k}
\sum\nolimits_{t=1}^{T_k^s}q_k^{t,s}
\leq\log_2\tfrac{\overline M}{M_0}.
\)
Counting one $F$ evaluation per completed inner iteration or rejected
inner trial, and one at $x_0$, gives the operator complexity. \hfill \Halmos
\end{proof}

Provided $\Lambda\geq\lambda^*+1$,
Proposition~\ref{thm:pfaos_oracle_complexity} shows that PFAOS requires at most
$\mathcal O(1/\sqrt{\varepsilon})$ constraint-gradient evaluations
and $\mathcal O(1/\varepsilon)$ operator evaluations to obtain an $\varepsilon$-approximate solution.
Two remarks clarify the parameter dependence of these bounds in the conclusion of Proposition \ref{thm:pfaos_oracle_complexity}.
First, the operator bound involves
$\overline M=2\max\{2\Lambda L,M\}$.
When $M\geq2\Lambda L$, the
$\overline{M}D_X^2/\varepsilon$ term becomes
$\mathcal O(MD_X^2/\varepsilon)$.
Note that such comparison between $M$ and $\Lambda L$ is not required for the preceding convergence or complexity results.
Second, the term $\overline M/\eta_0^{S_0}$ can be
large when $\eta_0^{S_0}$ is too small.
The condition $\eta_0^{S_0}\gtrsim\varepsilon/D_X^2$ allows this term
to be absorbed into $\mathcal O(\overline M D_X^2/\varepsilon)$.
Section~\ref{sec:special_first_iteration} modifies the search
at $k=1$ to remove this inverse dependence on $\eta_0^{S_0}$,
allowing arbitrarily small initial estimates.

\section{Multi-Constrained Parameter-Free Accelerated Operator Sliding (MC-PFAOS)}\label{sec:mc_pfaos}

Extending the single-constraint PFAOS approach directly to the setting with multiple constraints would require
solving each primal subproblem over $X$ intersected with $m$
halfspaces, which can be computationally expensive, especially for large $m$.
MC-PFAOS, presented in Algorithm~\ref{alg:mc-pfaos}, dualizes
these affine constraints and uses coupled primal--dual updates
requiring only projections onto $X$ and $\mathbb R_+^m$, respectively.

\begin{algorithm}[!htb]
\caption{Multi-constrained parameter-free accelerated operator sliding (MC-PFAOS)}
\label{alg:mc-pfaos}
\begin{algorithmic}[1]
\footnotesize
\State{\textbf{Initialization:} $\xub_0=\xlb_0^{S_0}=x_0\in X$, $\lamub_0=\lambda_0=0$, $L_0\in(0,L)$, \(\Lambda>0\),  and $M_0=\eta_0^{S_0}=2 \Lambda L_0$}
\For{$k=1,\dots,N$}
    \State{Set $s=1$ and $M_k^0 = M_{k-1}$}
    \State {Set \label{line:mcpfaos_s_loop}
        $\eta_k^s =  2^{s-2}\eta_{k-1}^{S_{k-1}}$
     and compute $\xlb_k^s = (1-\gamma_k^s)\xub_{k-1}+\gamma_k^sx_{k-1}$ with $\gamma_k^s$ defined by \eqref{eq:recursion_def}
    }
    \State{Set $t=1$,
$M_k^{0,s}=\max\{\eta_k^s,M_k^{s-1}\}$,
$T_k^{0,s}=M_k^{0,s}/\eta_k^s$,
$x_k^{0,s}=x_{k-1}$, 
$\lambda_k^{0,s}=\lambda_{k-1}$, 
$x_k^{-1,s}=x_{k-1}^{T_{k-1}^{S_{k-1}}-1,S_{k-1}}$
and
$\lambda_k^{-1,s}
=\lambda_{k-1}^{T_{k-1}^{S_{k-1}}-1,S_{k-1}}$
Here, 
$x_{0}^{T_{0}^{S_{0}}-1,S_{0}}=x_0$ and
$\lambda_{0}^{T_{0}^{S_{0}}-1,S_{0}}=0$.
}
    \State{Compute \label{line:mcpfaos_t_loop}
\begin{equation}\label{eq:mc_pfaos_subproblem}
    \begin{aligned}
	x_k^{t,s}=\arg\min_{x\in X}&
	\langle F(x_k^{t-1,s})+\theta_k^{t,s}[F(x_k^{t-1,s})-F(x_k^{t-2,s})],x\rangle +\tfrac{\eta_k^s}{2}\|x-x_{k-1}\|^2+\tfrac{p_k^{t,s}}{2}\|x-x_k^{t-1,s}\|^2\\
	&+\begin{cases}
	\langle \lambda_k^{t-1,s}+\theta_k^{t,s}(\lambda_k^{t-1,s}-\lambda_k^{t-2,s}), J(\xlb_k^s) x\rangle, \text{ if } t\ge2\\
	\langle \lambda_k^{0,s}, J(\xlb_k^s)x\rangle+\theta_k^{1,s}\langle \lambda_k^{0,s}-\lambda_k^{-1,s}, J(\xlb_{k-1}^{S_{k-1}})x\rangle \text{, if }t=1
	\end{cases}\\
	\lambda_k^{t,s} = \arg\min_{\lambda\in \mathbb{R}^m_+}& -\langle G(\xlb_k^s)+J(\xlb_k^s)(x_k^{t,s}-\xlb_k^s), \lambda\rangle+\tfrac{\eta_k^s}{2}\|\lambda-\lambda_{k-1}\|^2+\tfrac{p_k^{t,s}}{2}\|\lambda-\lambda_k^{t-1,s}\|^2
    \end{aligned}
    \end{equation}
    Here $\theta_1^{1,s}=1$, 
    $\theta_k^{1,s}=M_k^{1,s}/M_{k-1}$,
    $\forall k\ge2$, $\theta_k^{t,s}=M_k^{t,s}/M_k^{t-1,s}$,
    $\forall t\ge2$, and $p_k^{t,s} = 3M_k^{t,s}+\eta_k^sM_k^{t,s}\sum\nolimits_{i=1}^{t-1}(1/M_k^{i,s})$,  where $M_k^{t,s} = 2^{q_k^{t,s}}M_k^{t-1,s}$ and $q_k^{t,s}\ge0$ is the smallest integer such that 
    \(
    \|F(x_k^{t,s})-F(x_k^{t-1,s})\|\le M_k^{t,s}\|x_k^{t,s}-x_k^{t-1,s}\|
    \)  
    and 
    \(\|J(\xlb_k^s)^\top(\lambda_k^{t,s}-\lambda_k^{t-1,s})\|\le M_k^{t,s}\|\lambda_k^{t,s}-\lambda_k^{t-1,s}\|
    \)}
    \State{Set $T_k^{t,s}\gets t-1+2^{q_k^{t,s}}(T_k^{t-1,s}-t+1)$}
    \State{\textbf{if } $t<T_k^{t,s}$, \label{line:mcpfaos_line_8} \textbf{then } set $t\gets t+1$ and go to line \ref{line:mcpfaos_t_loop}; \textbf{else } set $T_k^s = T_k^{t,s}$, $M_k^s = M_k^{T_k^s,s}$,
    \begin{equation*}
            x_k^s = x_k^{T_k^s,s}
            ,\quad 
            \xhat_k^s = (\sum\nolimits_{i=1}^{T_k^s}1/M_k^{i,s})^{-1}\sum\nolimits_{i=1}^{T_k^s}x_k^{i,s}/M_k^{i,s}
            ,\quad
            \xub_k^s = (1-\gamma_k^s)\xub_{k-1}+\gamma_k^s\xhat_k^s,
            \text{ and }\lambda_k^s = \lambda_k^{T_k^s,s}
    \end{equation*}
    }
    \State{\textbf{if } \label{line:mcpfaos_middle_test}
    \(
    \|G(\xub_k^s)-G(\xlb_k^s)-J(\xlb_k^s)(\xub_k^s-\xlb_k^s)\|
    \le\tfrac{\eta_k^s}{2\Lambda\gamma_k^s}\|\xub_k^s-\xlb_k^s\|^2
    \)\text{ does not hold, }
    \textbf{then }    update $s\gets s+1$ and go to line \ref{line:mcpfaos_s_loop};  \textbf{else } set $S_k=s$, $M_k=M_k^{S_k}$,
    $x_k=x_k^{S_k}$, $\xhat_k = \xhat_k^{S_k}$, $\xub_k = \xub_k^{S_k}$,  $\lambda_k=\lambda_k^{S_k}$
    }
\EndFor
\State{\Return $\xub_N$}
\end{algorithmic}
\end{algorithm}

Algorithm~\ref{alg:mc-pfaos} retains the same outer-, middle-, and inner-loop structure as PFAOS. The updates of $\eta_k^s$, $\gamma_k^s$, $\theta_k^{t,s}$, 
and $T_k^{t,s}$, as well as the weighted averaging that defines
$\xhat_k^s$, follow the same mechanisms. 
We therefore focus on the modifications introduced in MC-PFAOS
relative to PFAOS.

First, at each middle trial $(k,s)$, MC-PFAOS linearizes the vector
constraint mapping $G(x)$ at $\xlb_k^s$ as
\(
G(\xlb_k^s)+J(\xlb_k^s)(x-\xlb_k^s),
\)
and reuses the same Jacobian $J(\xlb_k^s)$   throughout the   inner loop. Instead of imposing this affine model  as a primal
constraint, MC-PFAOS maintains a dual vector 
$\lambda_k^{t,s}\in\mathbb R_+^m$ for each inner iterate and performs coupled primal--dual hybrid-gradient-type proximal updates, resulting in strongly convex quadratic primal and dual problems over $X$ and $\mathbb R_+^m$, respectively.

Second, the primal update retains the operator extrapolation in PFAOS and adds dual extrapolation. For $t\geq2$, the extrapolated dual vector 
$\lambda_k^{t-1,s}+\theta_k^{t,s}(\lambda_k^{t-1,s}-\lambda_k^{t-2,s})$
is paired with the current Jacobian $J(\xlb_k^s)$. 
We treat the first inner iteration separately: the base vector $\lambda_k^{0,s}$ is paired
with $J(\xlb_k^s)$, whereas the inherited dual increment $\lambda_k^{0,s}-\lambda_k^{-1,s}$ is paired with
$J(\xlb_{k-1}^{S_{k-1}})$. This use of the previously accepted Jacobian allows the dual boundary terms to telescope across outer iterations.

Third, the inner line search adapts simultaneously to the local variation
of $F$ and the Jacobian--dual coupling. Unlike in PFAOS, where the local
estimate $M_k^{t,s}$ controls only the variation of $F$, in MC-PFAOS the
same estimate also controls the action of $J(\xlb_k^s)^\top$ on the
current dual increment. Accordingly, $M_k^{t,s}$ is increased until both
\(
\|F(x_k^{t,s})-F(x_k^{t-1,s})\|
\leq
M_k^{t,s}\|x_k^{t,s}-x_k^{t-1,s}\|
\)
and
\(
\|J(\xlb_k^s)^\top
(\lambda_k^{t,s}-\lambda_k^{t-1,s})\|
\leq
M_k^{t,s}
\|\lambda_k^{t,s}-\lambda_k^{t-1,s}\|
\)
hold. This is how the MC-PFAOS removes the prior knowledge of both the
Lipschitz constant and a uniform bound on $\|J(x)\|$ over $X$.

Finally, after the inner loop terminates, 
the middle-loop acceptance test is replaced by the vector Euclidean-norm
condition
\(
\|G(\xub_k^s)-G(\xlb_k^s)
-J(\xlb_k^s)(\xub_k^s-\xlb_k^s)\|
\leq
\tfrac{\eta_k^s}{2\Lambda\gamma_k^s}
\|\xub_k^s-\xlb_k^s\|^2.
\)
If this condition fails, $\eta_k^s$ is doubled and the next middle 
trial $(k,s+1)$ is performed; otherwise, the candidate is accepted as the output
of outer iteration $k$. 

\subsection{Convergence Analysis of MC-PFAOS}\label{subsec:mc_pfaos_convergence}
In this subsection, we establish the convergence properties of
Algorithm~\ref{alg:mc-pfaos}. We retain the definitions of
$\mathcal B_k^s$, $\Delta_k^{t,s}$, and $\Gamma_k^s$ from
Section~\ref{sec:pfaos}. For the dual iterates
$\lambda_k^{t,s}\in\mathbb R_+^m$ generated by
Algorithm~\ref{alg:mc-pfaos}, we define
$\hat\lambda_k^s$, $\hat\lambda_k$, and $\overline\lambda_k$
using the same averaging rules as in Section~\ref{sec:pfaos}.
We also denote the extrapolated dual variable by
$\lamtil_k^{t,s}:=\lambda_k^{t-1,s}
+\theta_k^{t,s}(\lambda_k^{t-1,s}-\lambda_k^{t-2,s})$.
In the multi-constraint setting, the affine constraint model and
primal--dual gap function become, for $x\in X$ and $\lambda\in\mathbb R_+^m$,
\[
\ell_k^s(z):=G(\xlb_k^s)+J(\xlb_k^s)(z-\xlb_k^s),
\qquad
\mathcal G_k(x,\lambda)
:=\langle F(x),\xub_k-x\rangle
+\langle G(\xub_k),\lambda\rangle
-\langle G(x),\overline\lambda_k\rangle.
\]

The horizon initialization and update are unchanged, so
Lemma~\ref{thm:t_loop_termination_invariance} also gives
$\eta_k^s\mathcal B_k^s=1$ for every middle trial in
Algorithm~\ref{alg:mc-pfaos}.
The convergence analysis follows similar inner-to-outer arguments as for
PFAOS, with additional terms arising from the dual updates.
Given this similarity, we defer the proofs of the results in
this and the next subsection to Appendix~\ref{app:mc_pfaos}
to streamline the presentation.

\begin{proposition}\label{thm:mc_pfaos_inner_analysis}
Fix any outer iteration $k\ge1$ and its middle trial $1\le s\le S_k$. For any $x\in X$ and vector $\lambda\in\mathbb R_+^m$, the iterates generated by Algorithm~\ref{alg:mc-pfaos} guarantee
\begin{equation}\label{eq:conclusion_mcpfaos_inner}
\mathcal{B}_k^s\langle F(x), \xhat_k^s-x\rangle 
    \le
    B_k^{s}
    + C_k^{s}
    + D_k^{s}
    -\sum\nolimits_{t=1}^{T_k^s}\tfrac{1}{M_k^{t,s}}\langle A_k^{t,s},x_k^{t,s}-x\rangle,
\end{equation}
where $A_k^{t,s}$ is the dual contribution to the primal objective:
$A_k^{1,s} = J(\xlb_k^s)^\top \lambda_k^{0,s} + \theta_k^{1,s}J(\xlb^{S_{k-1}}_{k-1})^\top(\lambda_k^{0,s}-\lambda_k^{-1,s})$
and 
$A_k^{t,s} = J(\xlb_k^s)^\top\lamtil_k^{t,s}$, for $t\ge 2$ and 
\begin{equation*}
\footnotesize
\begin{aligned}
       B_k^s :=&  \tfrac{1}{M_k^s}
    \langle\Delta_k^{T_k^s,s},x_k^s-x\rangle
        -\tfrac{\theta_k^{1,s}}{M_k^{1,s}}
        \langle\Delta_k^{0,s},x_{k-1}-x\rangle
        +\tfrac{p_k^{1,s}}{2M_k^{1,s}}
        (\|x-x_{k-1}\|^2+\|\lambda-\lambda_{k-1}\|^2)
        -\tfrac{\eta_k^s+p_k^{T_k^s,s}}{2M_k^s}
        (\|x-x_k^s\|^2+\|\lambda-\lambda_k^s\|^2) \\
        C_k^s:=& \tfrac{\theta_k^{1,s}}{M_k^{1,s}}\|\Delta_k^{0,s}\|\|x_k^{1,s}-x_{k-1}\| -\tfrac{p_k^{1,s}}{6M_k^{1,s}}\|x_k^{1,s}-x_{k-1}\|^2 - \tfrac{p_k^{T_k^s,s}}{6M_k^{T_k^s,s}}\|x_k^s-x_{k}^{T_k^s-1,s}\|^2 - \sum\nolimits_{t=1}^{T_k^s}\tfrac{p_k^{t,s}}{6M_k^{t,s}}\|x_k^{t,s}-x_k^{t-1,s}\|^2 \\
        D_k^s:=& \tfrac{1}{2}(\|x-x_{k-1}\|^2+\|\lambda-\lambda_{k-1}\|^2-\|\hat{x}_k^s-x_{k-1}\|^2) 
        + \sum\nolimits_{t=1}^{T_k^s}(\tfrac{1}{M_k^{t,s}}\langle\ell_k^s(x_k^{t,s}),\lambda_k^{t,s}-\lambda\rangle-\tfrac{p_k^{t,s}}{2M_k^{t,s}}\|\lambda_k^{t,s}-\lambda_k^{t-1,s}\|^2)
\end{aligned}
\end{equation*}
\end{proposition}

Combining the conclusion of Proposition~\ref{thm:mc_pfaos_inner_analysis} at $s=S_k$
with the middle-loop acceptance condition at line \ref{line:mcpfaos_middle_test} of Algorithm \ref{alg:mc-pfaos} yields the following
one-step bound on the outer primal--dual gap.

\begin{proposition}\label{thm:mc_pfaos_middle_gap}
For $k\ge1$, $x\in X$, and $\lambda\in\mathbb R_+^m$, the iterates generated by Algorithm~\ref{alg:mc-pfaos} satisfy
\begin{equation}\label{eq:mc_pfaos_outer_gap}
\mathcal G_k(x,\lambda)
-(1-\gamma_k^{S_k})\mathcal G_{k-1}(x,\lambda)
\le \tfrac{\gamma_k^{S_k}}{\mathcal{B}_k^{S_k}} (E_k+F_k+R'_k(\lambda)),
\end{equation}
where $R'_k(\lambda):=
\tfrac{1}{2}
 (\tfrac{\|\lambda\|}{\Lambda}-1 )
\|\xhat_k-x_{k-1}\|^2$ and 
\begin{equation*}
    \begin{aligned}
        E_k:=\;&
\tfrac{1}{M_k^{S_k}}
\langle\Delta_k^{T_k^{S_k},S_k},x_k-x\rangle
-\tfrac{\theta_k^{1,S_k}}{M_k^{1,S_k}}
\langle\Delta_k^{0,S_k},x_{k-1}-x\rangle
\\
&+\tfrac{1}{M_k^{S_k}}
\langle J(\xlb_k^{S_k})(x_k-x),\lambda_k-\lambda_k^{T_k^{S_k}-1,S_k}\rangle
-\tfrac{\theta_k^{1,S_k}}{M_k^{1,S_k}}
\langle J(\xlb_{k-1}^{S_{k-1}})(x_{k-1}-x),\lambda_{k-1}
-\lambda_{k-1}^{T_{k-1}^{S_{k-1}}-1,S_{k-1}}\rangle
\\
&+ 
\tfrac{1}{2}
(1+\tfrac{p_k^{1,S_k}}{M_k^{1,S_k}}
 )
(\|x-x_{k-1}\|^2+\|\lambda-\lambda_{k-1}\|^2)
-\tfrac{\eta_k^{S_k}+p_k^{T_k^{S_k},S_k}}
{2M_k^{S_k}}
(\|x-x_k\|^2+\|\lambda-\lambda_k\|^2),
\\
F_k:=\;&
\tfrac{\theta_k^{1,S_k}}{M_k^{1,S_k}}
\|\Delta_k^{0,S_k}\|\|x_k^{1,S_k}-x_{k-1}\|
+\tfrac{\theta_k^{1,S_k}}{M_k^{1,S_k}}
\|J(\xlb_{k-1}^{S_{k-1}})^\top (\lambda_{k-1}-\lambda_k^{-1,S_k})\|
\|x_k^{1,S_k}-x_{k-1}\|
\\
&-\tfrac{p_k^{1,S_k}}{3M_k^{1,S_k}}
\|x_k^{1,S_k}-x_{k-1}\|^2
-\tfrac{p_k^{T_k^{S_k},S_k}}{6M_k^{S_k}}
\|x_k-x_k^{T_k^{S_k}-1,S_k}\|^2
-\tfrac{p_k^{T_k^{S_k},S_k}}{2M_k^{S_k}}
\|\lambda_k-\lambda_k^{T_k^{S_k}-1,S_k}\|^2.
    \end{aligned}
\end{equation*}
\end{proposition}

With the help of Proposition~\ref{thm:mc_pfaos_middle_gap}, Theorem \ref{thm:mc_pfaos_global_convergence} next presents a global
primal--dual gap bound.

\begin{theorem}\label{thm:mc_pfaos_global_convergence}
For every $N\ge1$, $x\in X$, and $\lambda\in\mathbb R_+^m$,
the iterates generated by Algorithm~\ref{alg:mc-pfaos} satisfy
\begin{equation}\label{eq:mc_pfaos_global_gap}
\mathcal G_N(x,\lambda)
\le
\tfrac {\Gamma_N^{S_N}}{2}
 (
4\|x-x_0\|^2+4\|\lambda-\lambda_0\|^2
+ (\tfrac {\|\lambda\|}{\Lambda}-1 )
\sum\nolimits_{k=1}^N\|\xhat_k-x_{k-1}\|^2
 ).
\end{equation}
In particular, if $\Lambda\ge \|\lambda\|$, then
\(
\mathcal G_N(x,\lambda)
\le
2\Gamma_N^{S_N}
 (\|x-x_0\|^2+\|\lambda-\lambda_0\|^2 ).
\)
\end{theorem}

\subsection{Convergence Rate and Evaluation Complexity of MC-PFAOS}
\label{sec:mc_pfaos_complexity}
The preceding subsection provides an upper bound on the primal--dual gap
in terms of $\Gamma_N^{S_N}$. As in Section \ref{sec:pfaos_complexity}, we first bound $\Gamma_N^{S_N}$ and then
derive the evaluation complexities of MC-PFAOS.

\begin{lemma}\label{thm:mc_pfaos_Gamma_rate}
Algorithm~\ref{alg:mc-pfaos} guarantees that every middle-loop line
search terminates after finitely many trials and that 
\(
\eta_k^{S_k}/\gamma_k^{S_k}
\leq4\Lambda L
\), and 
\(
\Gamma_k^{S_k}
\leq
16\Lambda L/(k+1)^2
\)
,
\(
\forall k\geq1,
\)
\end{lemma}

Combining Lemma~\ref{thm:mc_pfaos_Gamma_rate} with
Theorem~\ref{thm:mc_pfaos_global_convergence} yields an explicit
primal--dual gap rate. We next translate it into bounds on the weak VI gap
and constraint violation in Corollary \ref{thm:mc_pfaos_duality_gap_rate}.

\begin{corollary}\label{thm:mc_pfaos_duality_gap_rate}
Let $(x^*,\lambda^*)$ satisfy
\eqref{eq:standing_primal_dual_pair}. For every $N\geq1$, we have
\begin{equation*}
\sup_{x\in\widetilde X}
\langle F(x),\xub_N-x\rangle
\leq2\Gamma_N^{S_N}D_X^2
\leq\tfrac{32\Lambda LD_X^2}{(N+1)^2},\text{ and }
(\Lambda-\|\lambda^*\|)\|[G(\xub_N)] 
_+\|
\leq2\Gamma_N^{S_N}(D_X^2+\Lambda^2)
\leq\tfrac{32\Lambda L(D_X^2+\Lambda^2)}{(N+1)^2}.
\end{equation*}
In particular, if $\Lambda\geq\|\lambda^*\|+1$, then $\xub_N$
is an $\varepsilon$-approximate solution of
\eqref{eq:problem_of_interest} whenever
\(
2\Gamma_N^{S_N}(D_X^2+\Lambda^2)\leq\varepsilon,
\)
which is guaranteed by
\(
N+1\geq
\sqrt{\tfrac{32\Lambda L(D_X^2+\Lambda^2)}{\varepsilon}}.
\)
\end{corollary}

Corollary~\ref{thm:mc_pfaos_duality_gap_rate} gives a sufficient number of outer
iterations to obtain an $\varepsilon$-approximate solution. Using this result, Proposition \ref{thm:mc_pfaos_oracle_complexity} counts the number of  Jacobian and operator evaluations.

\begin{proposition}\label{thm:mc_pfaos_oracle_complexity}
Suppose $\Lambda\geq\|\lambda^*\|+1$. For any $\varepsilon>0$,
Algorithm~\ref{alg:mc-pfaos} obtains an $\varepsilon$-approximate
solution of \eqref{eq:problem_of_interest} using at most
\(
3+2\sqrt{\tfrac{32\Lambda L(D_X^2+\Lambda^2)}{\varepsilon}}
+\lceil\log_2\tfrac{2L}{L_0}\rceil
\)
evaluations of $J$ and
\(
1+\tfrac{4\overline M_{\rm MC}}{\eta_0^{S_0}}
+\tfrac{8\overline M_{\rm MC}(D_X^2+\Lambda^2)}{\varepsilon}
+\lceil\log_2\tfrac{\overline M_{\rm MC}}{M_0}\rceil
\)
evaluations of $F$, where
\(
\overline M_{\rm MC}:=2\max\{2\Lambda L,M,\max_{x\in X}\|J(x)\|\}.
\)
These bounds include evaluations at initialization and in
accepted and rejected line-search trials.
\end{proposition}

Provided $\Lambda\geq\|\lambda^*\|+1$,
Proposition~\ref{thm:mc_pfaos_oracle_complexity} shows that 
MC-PFAOS requires at most
$\mathcal O(1/\sqrt{\varepsilon})$ Jacobian evaluations
and $\mathcal O(1/\varepsilon)$ operator evaluations for obtaining an $\varepsilon$-approximate solution. 
Three remarks apply here.
First, no prior bound on $\|J(x)\|$ is required by the algorithm.
Second, when
$M\geq\max\{2\Lambda L,\max_{x\in X}\|J(x)\|\}$,
the $\overline M_{\rm MC}(D_X^2+\Lambda^2)/\varepsilon$
term in the operator bound becomes
$\mathcal O(M(D_X^2+\Lambda^2)/\varepsilon)$.
This condition is not required for the preceding convergence
or complexity results.
Third, the term $\overline M_{\rm MC}/\eta_0^{S_0}$
is absorbed into $\mathcal O(\overline{M}_{MC}(D_X^2+\Lambda^2)/\varepsilon)$  when
$\eta_0^{S_0}\gtrsim\varepsilon/(D_X^2+\Lambda^2)$.
 Such
initialization restriction of $\eta_0^{S_0}$ is removed in Section~\ref{sec:special_first_iteration} through a modified search at $k=1$.

\section{Modified First-Iteration Parameter Search}
\label{sec:special_first_iteration}

The operator-complexity bounds in
Propositions~\ref{thm:pfaos_oracle_complexity} and
\ref{thm:mc_pfaos_oracle_complexity} contain the 
terms $\overline M/\eta_0^{S_0}$ and $\overline M_{\rm MC}/\eta_0^{S_0}$, respectively,
which can be arbitrarily large when the initial estimate $\eta_0^{S_0}$ is small.
Absorbing these terms directly into the desired $\mathcal O(1/\varepsilon)$ operator complexity bounds
requires, respectively, $\eta_0^{S_0}\gtrsim\varepsilon/D_X^2$ and $\eta_0^{S_0}\gtrsim\varepsilon/(D_X^2+\Lambda^2).$
However, these conditions are not guaranteed by the initialization $\eta_0^{S_0}=2\Lambda L_0$ with $0<L_0\leq L$.
This section modifies the search order for $\eta_1^s$ at $k=1$ in
Algorithms~\ref{alg:pfaos} and~\ref{alg:mc-pfaos}
to achieve the desired $\mathcal O(1/\varepsilon)$ operator
complexity without these initialization requirements.
We retain the names PFAOS and MC-PFAOS because it modifies only the search order of the first outer iteration. 

\subsection{First-Iteration Modifications}

Algorithms~\ref{alg:special_k1_pfaos} and
\ref{alg:special_k1_mcpfaos} replace the $k=1$ iteration of
Algorithms \ref{alg:pfaos} and \ref{alg:mc-pfaos}, respectively.
At $s=1$, an upward search doubles the operator and curvature
estimates together from $M_0=\eta_0^{S_0}$ until all inner and middle
acceptance conditions hold. 
For $s\geq2$, a downward search halves $\eta_1^s$ between
middle trials while retaining the original inner-loop rules.
It returns $\xub_1^s$ if the $s$-th trial is accepted and satisfies
the accuracy certificates of \(
\tfrac32\eta_1^sD_X^2\leq\varepsilon
\) or  \(
2\eta_1^s(D_X^2+\Lambda^2)\le \varepsilon\).
If a downward trial fails the middle test, then its preceding
trial is retained and the original algorithm resumes at $k=2$. 
Here, each upward or downward trial starts from
$x_1^{0,s}=x_1^{-1,s}=x_0$, with
$\lambda_1^{0,s}=\lambda_1^{-1,s}=0$ for MC-PFAOS.
Since $\gamma_1^s=1$ and $\xlb_1^s=x_0$ throughout both searches,
one evaluation of $\nabla g(x_0)$ or $J(x_0)$ suffices.

\begin{algorithm}[!htb]
\caption{Modified first outer iteration for PFAOS}
\label{alg:special_k1_pfaos}
\footnotesize
\begin{algorithmic}[1]
\State{Same initialization and input as in Algorithm \ref{alg:pfaos}}
\State{\textbf{Upward search:} Set $k=1$, $s=1$; Run lines \ref{line:pfaos_s_loop} - \ref{alg1_line:8} in Algorithm \ref{alg:pfaos} with fixed $T_1^{1,1}=1$ and $\eta_1^1=M_1^{1,1}=2^{q}M_0$, where $q\ge0$ is the smallest integer such that both  
\(
\|F(x_1^{1,1})-F(x_0)\|\leq M_1^{1,1}\|x_1^{1,1}-x_0\|
\)
and 
\(
g(\xub_1^1)-g(x_0)-\langle\nabla g(x_0),\xub_1^1-x_0\rangle\leq\tfrac{\eta_1^1}{2\Lambda}\|\xub_1^1-x_0\|^2
\) hold
}
\State{\textbf{If }
\(
\tfrac32\eta_1^1D_X^2\leq\varepsilon
\)
\textbf{, then } set $S_1=1$, $M_1=M_1^{S_1}$, $x_1=x_1^{S_1}$,
$\xhat_1=\xhat_1^{S_1}$ and output
$\xub_1=\xub_1^{S_1}$; \textbf{ Else } set $s=2$}
\State{\textbf{Downward search: }\label{line:special_k1_downsearch}
Run line \ref{line:pfaos_s_loop} - \ref{alg1_line:8} in Algorithm \ref{alg:pfaos} with $\eta_1^s=\eta_1^{s-1}/2$}

\State \textbf{If }  \text{the middle test condition }
$g(\xub_1^s)\le g(x_0)+\langle\nabla g(x_0), \xub_1^s-x_0\rangle+\tfrac{\eta_1^s}{2\Lambda}\|\xub_1^s-x_0\|^2$\; \text{does not hold,} 
\textbf{ then } Set $S_1=s-1$, 
$M_1=M_1^{S_1}$, $x_1=x_1^{S_1}$,
$\xhat_1=\xhat_1^{S_1}$, 
$\xub_1=\xub_1^{S_1}$ and resume Algorithm \ref{alg:pfaos} at $k=2$;
 \textbf{Else if } \(
\tfrac32\eta_1^sD_X^2\leq\varepsilon
\) \textbf{ then } set $S_1=s$ 
$M_1=M_1^{S_1}$, $x_1=x_1^{S_1}$,
$\xhat_1=\xhat_1^{S_1}$ and output 
$\xub_1=\xub_1^{S_1}$;  \textbf{Else }
set $s\gets s+1$ and go to line \ref{line:special_k1_downsearch} 
\end{algorithmic}
\end{algorithm}

\begin{algorithm}[!htb]
\caption{Modified first outer iteration for MC-PFAOS}
\label{alg:special_k1_mcpfaos}
\footnotesize
\begin{algorithmic}[1]
\State{Same initialization and input as in Algorithm \ref{alg:mc-pfaos}}
\State{\textbf{Upward search:}\label{line:special_k1_upward}
Set $k=1$, $s=1$; Run lines \ref{line:mcpfaos_s_loop}-\ref{line:mcpfaos_line_8} in Algorithm \ref{alg:mc-pfaos} with fixed $T_1^{1,1}=1$ and $\eta_1^1=M_1^{1,1}=2^qM_0$, where $q\ge0$ is the smallest integer such that 
\(
\|F(x_1^{1,1})-F(x_0)\|
\leq M_1^{1,1}\|x_1^{1,1}-x_0\|
\),
\(
\|J(x_0)^\top(\lambda_1^{1,1}-\lambda_0)\|
\leq M_1^{1,1}\|\lambda_1^{1,1}-\lambda_0\|
\) and 
\(
\|G(\xub_1^1)-G(x_0)-J(x_0)(\xub_1^1-x_0)\|\le \tfrac{\eta_1^1}{2\Lambda}\|\xub_1^1-x_0\|^2
\) hold}
\State{\textbf{If } \(
2\eta_1^1(D_X^2+\Lambda^2)\le \varepsilon\) \textbf{, then } set $S_1=1$, $M_1=M_1^{S_1}$,  $\lambda_1=\lambda_1^{S_1}$,
$x_1=x_1^{S_1}$,
$\xhat_1=\xhat_1^{S_1}$ and output
$\xub_1=\xub_1^{S_1}$; \textbf{ Else } set $s=2$}
\State{\textbf{Downward search:}\label{line:special_k1_downward_mc}
Run line \ref{line:mcpfaos_s_loop}-\ref{line:mcpfaos_line_8} in Algorithm \ref{alg:mc-pfaos} with $\eta_1^s=\eta_1^{s-1}/2$}
\State{\textbf{If } \text{the middle test condition } \(
\|G(\xub_1^s)-G(x_0)-J(x_0)(\xub_1^s-x_0)\|\le \tfrac{\eta_1^s}{2\Lambda}\|\xub_1^s-x_0\|^2
\) \text{ does not hold, }\textbf{then } set $S_1=s-1$, $M_1=M_1^{S_1}$,  $\lambda_1=\lambda_1^{S_1}$,
$x_1=x_1^{S_1}$,
$\xhat_1=\xhat_1^{S_1}$
$\xub_1=\xub_1^{S_1}$ and resume Algorithm \ref{alg:mc-pfaos} at $k=2$; 
\textbf{Else if}  \(
2\eta_1^s(D_X^2+\Lambda^2)\le \varepsilon\) \textbf{, then } set $S_1=s$, $M_1=M_1^{S_1}$,  $\lambda_1=\lambda_1^{S_1}$,
$x_1=x_1^{S_1}$,
$\xhat_1=\xhat_1^{S_1}$ and output
$\xub_1=\xub_1^{S_1}$; \textbf{ Else } set $s\gets s+1$ and go to line \ref{line:special_k1_downward_mc}}
\end{algorithmic}
\end{algorithm}
\subsection{Termination and Convergence Guarantees}
We first establish the finite termination of Algorithms \ref{alg:special_k1_pfaos} and \ref{alg:special_k1_mcpfaos}, and estimate the evaluation bounds of the modified first outer iteration in Lemma \ref{thm:special_k1_cost}.
\begin{lemma}\label{thm:special_k1_cost}
Algorithms~\ref{alg:special_k1_pfaos} and \ref{alg:special_k1_mcpfaos} terminate after finitely many trials.
If it continues to $k=2$, then $\eta_1^{S_1}<2\Lambda L$.
The first outer iteration requires one evaluation of $\nabla g$ or $J$
and at most
\(
2+\tfrac{6\overline M D_X^2}{\varepsilon}
+\lceil\log_2\tfrac{\overline M}{M_0}\rceil
\)
or
\(
2+\tfrac{8\overline M_{\rm MC}(D_X^2+\Lambda^2)}{\varepsilon}
+\lceil\log_2\tfrac{\overline M_{\rm MC}}{M_0}\rceil
\)
evaluations of $F$ for PFAOS and MC-PFAOS, respectively.
Both bounds include evaluations at initialization and in
accepted and rejected  trials.
\end{lemma}

\begin{proof}{Proof: }
By the sufficient acceptance conditions and the doubling rule at
$s=1$, the upward search terminates with a common estimate at most
$\overline M$ or $\overline M_{\rm MC}$, respectively.
For $s\geq2$, halving $\eta_1^s$ preserves
$M_1^{t,s}\leq\overline M$ for PFAOS and
$M_1^{t,s}\leq\overline M_{\rm MC}$ for MC-PFAOS
by the arguments in Propositions~\ref{thm:pfaos_oracle_complexity}
and~\ref{thm:mc_pfaos_oracle_complexity},
and ensures $S_1<\infty$, since each inner loop is finite and
either a middle test fails or an accepted trial eventually
satisfies the accuracy certificate. 
If the method continues to $k=2$, the final middle trial $s$ failed,
so $\eta_1^s<\Lambda L$.
Retaining $S_1=s-1$ therefore gives
$\eta_1^{S_1}=2\eta_1^s<2\Lambda L$ for both methods.
Since $\xlb_1^s=x_0$ throughout the search, one evaluation of
$\nabla g(x_0)$ or $J(x_0)$ suffices.
For the PFAOS operator count, let $s$ be the final tested middle-trial
index, including the rejected downward trial, if any.
If $s\geq2$, the preceding accepted trial did not satisfy the accuracy
certificate, so
\(
\eta_1^s=\eta_1^{s-1}/2>\varepsilon/(3D_X^2).
\)
Using $\eta_1^i\mathcal B_1^i=1$ and $M_1^{t,i}\leq\overline M$ gives
\(
\sum\nolimits_{i=1}^{s}T_1^i
\leq\overline M\sum\nolimits_{i=1}^{s}\tfrac1{\eta_1^i}
<\tfrac{2\overline M}{\eta_1^s}
<\tfrac{6\overline M D_X^2}{\varepsilon}.
\)
Here, $T_1^1=1$ counts only the accepted upward attempt.
For $s=1$, the sum equals one, so it is at most
$1+6\overline M D_X^2/\varepsilon$ in either case.
Each rejected upward attempt or inner backtracking trial doubles
the nondecreasing operator estimate, so their total number is at most
$\lceil\log_2(\overline M/M_0)\rceil$.
Adding these evaluations and $F(x_0)$ to the completed-iteration
count proves the PFAOS bound.
For MC-PFAOS, the same argument uses $\overline M_{\rm MC}$ and
$\eta_1^s>\varepsilon/[4(D_X^2+\Lambda^2)]$ for $s\geq2$,
yielding the stated bound. \hfill \Halmos
\end{proof}

Next, we establish in Lemma \ref{thm:special_k1_convergence} the gap and parameter bounds for the modified methods.
\begin{lemma}\label{thm:special_k1_convergence}
The bounds in
Theorems~\ref{thm:pfaos_global_convergence} and
\ref{thm:mc_pfaos_global_convergence} remain valid with the
first-iteration modification in Algorithms~\ref{alg:special_k1_pfaos} and \ref{alg:special_k1_mcpfaos}.
If $\Lambda\ge \|\lambda^*\|+1$, then any point returned at $k=1$ is an $\varepsilon$-approximate solution
of \eqref{eq:problem_of_interest}.
If the method continues to $k=2$, then
\(
\eta_k^{S_k}/\gamma_k^{S_k}\leq4\Lambda L
\)
and
\(
\Gamma_k^{S_k}\leq16\Lambda L/(k+1)^2
\)
for $k\geq1$.
\end{lemma}
\begin{proof}{Proof: }
Each retained trial satisfies the original inner and middle acceptance conditions and uses the original initialization, inner updates, and
averaging rules.
The invariant $\eta_1^{S_1}\mathcal B_1^{S_1}=1$ holds because the accepted upward attempt has $M_1^{1,1}=\eta_1^1$ and $T_1^1=1$,
while downward trials use the original horizon update.
So, the arguments in Theorems~\ref{thm:pfaos_global_convergence} and
\ref{thm:mc_pfaos_global_convergence} remain valid.
Since $\Gamma_1^{S_1}=\eta_1^{S_1}$, the stopping tests at $k=1$ ensure the claimed accuracy by Corollaries~\ref{thm:pfaos_duality_gap_rate} and \ref{thm:mc_pfaos_duality_gap_rate}.
If the method continues to $k=2$, Lemma~\ref{thm:special_k1_cost}
gives $\eta_1^{S_1}<2\Lambda L$.
Together with $\gamma_1^{S_1}=1$, this provides the base case.
Then induction and telescoping arguments in Lemmas~\ref{thm:pfaos_Gamma_rate} and \ref{thm:mc_pfaos_Gamma_rate} yield the stated bounds. \hfill \Halmos
\end{proof}

\subsection{Overall Evaluation Complexity}
Combining the preceding Lemmas \ref{thm:special_k1_cost} and \ref{thm:special_k1_convergence} with the counting arguments in
Propositions~\ref{thm:pfaos_oracle_complexity} and
\ref{thm:mc_pfaos_oracle_complexity} yields the following overall
evaluation bounds for the modified methods.
\begin{theorem}\label{thm:special_k1_complexity}
Suppose 
$\Lambda\geq\|\lambda^*\|+1$.
For any $\varepsilon>0$, both methods with the modification in
Algorithm~\ref{alg:special_k1_pfaos} and \ref{alg:special_k1_mcpfaos} obtain an $\varepsilon$-approximate
solution of \eqref{eq:problem_of_interest}.
PFAOS requires at most
\(
2+2\sqrt{\tfrac{24\Lambda L D_X^2}{\varepsilon}}
+\lceil\log_2(1+\tfrac{6\Lambda L D_X^2}{\varepsilon})\rceil
\)
evaluations of $\nabla g$ and
\(
2+\tfrac{12\overline M D_X^2}{\varepsilon}
+2\lceil\log_2\tfrac{\overline M}{M_0}\rceil
\)
evaluations of $F$.
For MC-PFAOS, the corresponding bounds are
\(
2+2\sqrt{\tfrac{32\Lambda L(D_X^2+\Lambda^2)}{\varepsilon}}
+\lceil\log_2(1+\tfrac{8\Lambda L(D_X^2+\Lambda^2)}{\varepsilon})\rceil
\)
evaluations of $J$ and
\(
2+\tfrac{16\overline M_{\rm MC}(D_X^2+\Lambda^2)}{\varepsilon}
+2\lceil\log_2\tfrac{\overline M_{\rm MC}}{M_0}\rceil
\)
evaluations of $F$.
These bounds include evaluations at initialization and in
accepted and rejected line-search trials.
\end{theorem}

\begin{proof}{Proof: }
We first consider PFAOS. Let $N_\varepsilon\geq1$ be the first outer
iteration satisfying
\(
\tfrac32\Gamma_{N_\varepsilon}^{S_{N_\varepsilon}}D_X^2
\leq\varepsilon.
\)
By Lemma~\ref{thm:special_k1_convergence} and
Corollary~\ref{thm:pfaos_duality_gap_rate},
$\xub_{N_\varepsilon}$ is an $\varepsilon$-approximate solution.
If $N_\varepsilon=1$, the evaluation bounds follow from
Lemma~\ref{thm:special_k1_cost}, allowing one additional initialization
evaluation of $\nabla g$.
Hence, assume $N_\varepsilon\geq2$.
Lemma~\ref{thm:special_k1_convergence} then gives
\(
N_\varepsilon
\leq1+\sqrt{\tfrac{24\Lambda L D_X^2}{\varepsilon}}.
\)
For $\nabla g$ evaluations, the retained trial at $k=1$
did not satisfy the accuracy certificate, so
$\eta_1^{S_1}>2\varepsilon/(3D_X^2)$.
Together with $\eta_{N_\varepsilon}^{S_{N_\varepsilon}}\leq4\Lambda L$
,
the $\eta_k^s$ updates give
\begin{equation*}
\sum\nolimits_{k=2}^{N_\varepsilon}S_k
=2(N_\varepsilon-1)
+\log_2\tfrac{\eta_{N_\varepsilon}^{S_{N_\varepsilon}}}{\eta_1^{S_1}}
<2(N_\varepsilon-1)
+\log_2(1+\tfrac{6\Lambda L D_X^2}{\varepsilon}).
\end{equation*}
Each middle trial requires one gradient evaluation.
Adding the evaluation at $x_0$, as in Proposition~\ref{thm:pfaos_oracle_complexity},
and using the bound on $N_\varepsilon$ yields the stated
constraint-gradient complexity.
For the operator evaluations, the invariant and the geometric-series
argument in Proposition~\ref{thm:pfaos_oracle_complexity} give
\begin{equation*}
\sum\nolimits_{k=2}^{N_\varepsilon}
\sum\nolimits_{s=1}^{S_k}T_k^s
<4\overline M\sum\nolimits_{k=1}^{N_\varepsilon-1}
\tfrac1{\eta_k^{S_k}}
=\tfrac{4\overline M}
{\Gamma_{N_\varepsilon-1}^{S_{N_\varepsilon-1}}}
<\tfrac{6\overline M D_X^2}{\varepsilon}.
\end{equation*}
Retaining the preceding trial at $k=1$ may lower the operator estimate
to $M_1=M_1^{S_1}$.
For $k\ge2$, the estimates are nondecreasing from
$M_1\geq M_0$ and remain bounded by $\overline M$.
So, there are at most
$\lceil\log_2(\overline M/M_0)\rceil$
additional rejected inner trials.
Combining these counts with Lemma~\ref{thm:special_k1_cost}
gives the stated operator complexity.
For MC-PFAOS, the same argument replaces $\overline M$ and $D_X^2$
by $\overline M_{\rm MC}$ and $D_X^2+\Lambda^2$, respectively,
and uses the stopping condition
\(
2\Gamma_{N_\varepsilon}^{S_{N_\varepsilon}}
(D_X^2+\Lambda^2)\leq\varepsilon.
\)
\hfill
\Halmos
\end{proof}

Theorem~\ref{thm:special_k1_complexity} indicates that the special
first-iteration searches in
Algorithms~\ref{alg:special_k1_pfaos} and
\ref{alg:special_k1_mcpfaos} allow arbitrarily small positive
initial estimates $\eta_0^{S_0}$ while preserving
$\mathcal O(1/\sqrt{\varepsilon})$ constraint-gradient or
Jacobian complexity and $\mathcal O(1/\varepsilon)$ operator
complexity for both PFAOS and MC-PFAOS.
These guarantees use accuracy tests involving $D_X$, both during
the special first iteration and when determining the number
of outer iterations $N$.
Since $\gamma_1^s=1$ gives $\Gamma_1^s=\eta_1^s$,
both types of tests can be written as
\(
\tfrac32\Gamma_N^{S_N}D_X^2\leq\varepsilon
\)
for PFAOS and 
\(
2\Gamma_N^{S_N}(D_X^2+\Lambda^2)\leq\varepsilon\) for
MC-PFAOS, respectively, with $N=1$ and $S_1=s$
when testing an accepted first-iteration trial.
However, $D_X$ is not available.
Section~\ref{sec:stopping_criteria} therefore replaces these tests
with an automatic stopping rule that certifies an
$\varepsilon$-approximate solution without requiring knowledge of $D_X$
and preserves these evaluation-complexity bounds.
\section{Computable Stopping Criteria}\label{sec:stopping_criteria}
In this section, we construct a computable weak VI gap certificate from quantities
generated by PFAOS and MC-PFAOS.
Combined with the observed constraint violation, it provides
an automatic stopping rule that requires no knowledge of $D_X$
and preserves the established complexity bounds. The proofs of the results in this section are provided in
Appendix~\ref{app:stopping_criteria}.
\subsection{A Computable Gap Certificate}
We use vector notations for both methods, with $m=1$ for PFAOS, and
define 
\begin{equation}\label{eq:stopping_affine_model}
\mathcal C_N:={} [\max_{x\in X}\Phi_N(x)  ]_+ \text{ with }\Phi_N(x):={}\Gamma_N^{S_N}
\sum\nolimits_{k=1}^N\sum\nolimits_{t=1}^{T_k^{S_k}}
[
\langle F(x_k^{t,S_k}),x_k^{t,S_k}-x\rangle
-\langle\lambda_k^{t,S_k},\ell_k^{S_k}(x)\rangle
]/M_k^{t,S_k}.
\end{equation}
Here, evaluating $\mathcal C_N$ requires only 
linear optimization over the simple set $X$. It uses stored quantities and requires no additional evaluations of $F$, $\nabla g$ or $J$.
The following proposition shows that this certificate bounds
the weak VI gap and is controlled by the preceding
$D_X$-based bounds.
\begin{proposition}\label{thm:stopping_certificate}
For both methods, with or without the special first-iteration search,
every retained outer iterate satisfies
\(
\sup_{x\in\widetilde X}\langle F(x),\xub_N-x\rangle
\leq\mathcal C_N.
\)
Moreover,
\(
\mathcal C_N\leq\tfrac32\Gamma_N^{S_N}D_X^2
\)
for PFAOS and
\(
\mathcal C_N\leq2\Gamma_N^{S_N}D_X^2
\)
for MC-PFAOS.
\end{proposition}

\subsection{Stopping Rule and Complexity }

With the help of
Proposition~\ref{thm:stopping_certificate},
we can replace the preceding $D_X$-based accuracy
tests by the following rule: after obtaining an accepted outer update,
return $\xub_N$ if
\begin{equation}\label{eq:computable_stopping_rule}
\mathcal C_N\leq\varepsilon,
\qquad
\|[G(\xub_N)]_+\|\leq\varepsilon.
\end{equation}
For the special first iteration in
Algorithms~\ref{alg:special_k1_pfaos} and
\ref{alg:special_k1_mcpfaos}, replace each accuracy test
with \eqref{eq:computable_stopping_rule}.
For a trial $s$ that passes the inner and middle acceptance
conditions, evaluate the certificate with $N=1$, $S_1=s$,
and $\Gamma_1^s=\eta_1^s$, using only that trial's inner iterates.
All other updates remain unchanged.
The following proposition establishes the accuracy and
complexity guarantees of this replacement.

\begin{proposition}\label{thm:stopping_complexity}
For any $\Lambda>0$, any point returned by
\eqref{eq:computable_stopping_rule} is an $\varepsilon$-approximate
solution of \eqref{eq:problem_of_interest}.
If $\Lambda\geq\|\lambda^*\|+1$, \eqref{eq:computable_stopping_rule} preserves the
respective complexities in
Propositions~\ref{thm:pfaos_oracle_complexity} and
\ref{thm:mc_pfaos_oracle_complexity}, and
Theorem~\ref{thm:special_k1_complexity}.
\end{proposition}

Therefore, the stopping rules in \eqref{eq:computable_stopping_rule} certify an $\varepsilon$-approximate
solution without requiring $D_X$, while retaining the
established complexity bounds.

\section{Numerical Experiments}\label{sec:numerical_experiments}
To our knowledge, no established parameter-free method directly applies
to CVI setting, and adapting parameter-dependent CVI methods is
nontrivial. We therefore benchmark PFAOS and MC-PFAOS against
unconstrained VI methods using projections onto the full feasible set.
To assess robustness, we generate  instances across multiple
random seeds, ensuring the required monotonicity and feasibility
conditions as described below.

\subsection{Experiment Setting}
\subsubsection{Test Instances}
\label{subsec:numerical_instances}
We generate randomized instances inspired by bilinear games with
a singular positive semidefinite (PSD) regularizer.
Let $X=\{x\in\mathbb R^n:\|x\|_2\leq 2\}$
and generate operator 
$F(x)=(R-R^\top+H)x-d$ and quadratic constraints
$g_i(x)=\tfrac12x^\top A_i x+b_i^\top x+c_i$,
where $R,H,A_i\in\mathbb R^{n\times n}$,
$d,b_i\in\mathbb R^n$, and $c_i\in\mathbb R$.
Uniform draws are used to generate $R,d,b_i$ and the positive
eigenvalues of $H$ and $A_i$, with
$\lfloor n/2\rfloor$ eigenvalues fixed at zero.
We construct $H$ as a diagonal PSD matrix and  $A_i$ as symmetric
PSD matrices
by conjugating the corresponding diagonal matrix with a random
orthogonal matrix.
Thus, $F$ is nonconservative and monotone but not strongly monotone, and the constraints
are convex but not strongly convex.
The constants $c_i$ are chosen so that a common point
$x^{\mathrm sl}\in\operatorname{ri}(X)$ satisfies
$g_i(x^{\mathrm sl})<0$ for all $i=1,\ldots,m$.
Consequently, Slater's condition in
Assumption~\ref{ass:pridual_slater} holds, and the full feasible set
$\widetilde X=X\cap\{x:g_i(x)\leq0,\ i=1,\ldots,m\}$
is nonempty, compact, and convex, ensuring that the VI admits a solution.
The instance-generation and implementation code is available at \url{https://github.com/qluorocky/CVI_IJOC}.
\subsubsection{Benchmark Methods}
\label{subsec:numerical_methods}
We adapt two advanced first-order algorithms that achieve competitive sample complexity for unconstrained VIs. 
\begin{enumerate}
    \item Benchmark 1 applies the operator-extrapolation (OE) method in \citet{kotsalis2020simplea}.  
    For a fixed  estimate $\widehat M>0$, OE updates 
    \(x_{k+1} =\Pi_{\widetilde X} ( x_k- (2F(x_k)-F(x_{k-1}))/(2\widehat M) )\) at iteration $k$.
    For each  instance, we run fixed-$\widehat M$ trajectories in parallel with trial values spanning $[0.1M,10M]$. 
    
    \item Benchmark 2  applies the golden-ratio algorithm in  \citet{malitsky2020golden}, which is a parameter-free method. 
    Given 
   initializations, the method updates 
        \(
        \bar x_k  =((\phi-1)x_k+\bar x_{k-1})/\phi 
        \)
        and   
        \(
        x_{k+1} =\Pi_{\widetilde X}(\bar x_k-\lambda_k F(x_k))
        \), where $\varphi=(1+\sqrt{5})/2$,  $\phi\in(1,\varphi]$ and \(
        \lambda_k =\min \{ \  (\tfrac{1}{\phi}+\tfrac{1}{\phi^{2}})\lambda_{k-1}  ,\ \tfrac{\phi^2  }{4\lambda_{k-2}} \tfrac{\|x_k-x_{k-1}\|^2} {\|F(x_k)-F(x_{k-1})\|^2},\ \bar\lambda  \},
        \).
\end{enumerate}

Both benchmarks compute projections  $\Pi_{\widetilde X} $ with a generic convex solver at every iteration. 
\subsubsection{Evaluation Protocol}
\label{subsec:numerical_protocol}
All methods run with the same set of realized instances and initial points.
The runs were performed on an Apple M4 MacBook Pro (10 CPU cores and 16 GB memory), using  CVXPY 1.9.3 and Clarabel 0.11.1.
For a reported average $\bar x$, we report both the VI gap  $ \operatorname{Gap}(\bar x) = [\sup_{y\in\widetilde X}
 \langle F(y),\bar x-y\rangle ]_+,$ and constraint violation  $ V(\bar x)=\|[G(\bar x)]_+\|_2. $
 We denote the combined error $ \widehat E =    \max\{ \operatorname{Gap}(\bar{x}), \|[G(\bar{x})]_+\|_2\}$
The common stopping target is $10^{-5}$ and we also report their runtime along the convergence paths. 

\subsection{Results}\label{subsec:numerical_results}
\begin{figure}[!htbp]
  \centering
  \includegraphics[width=0.85 \linewidth]{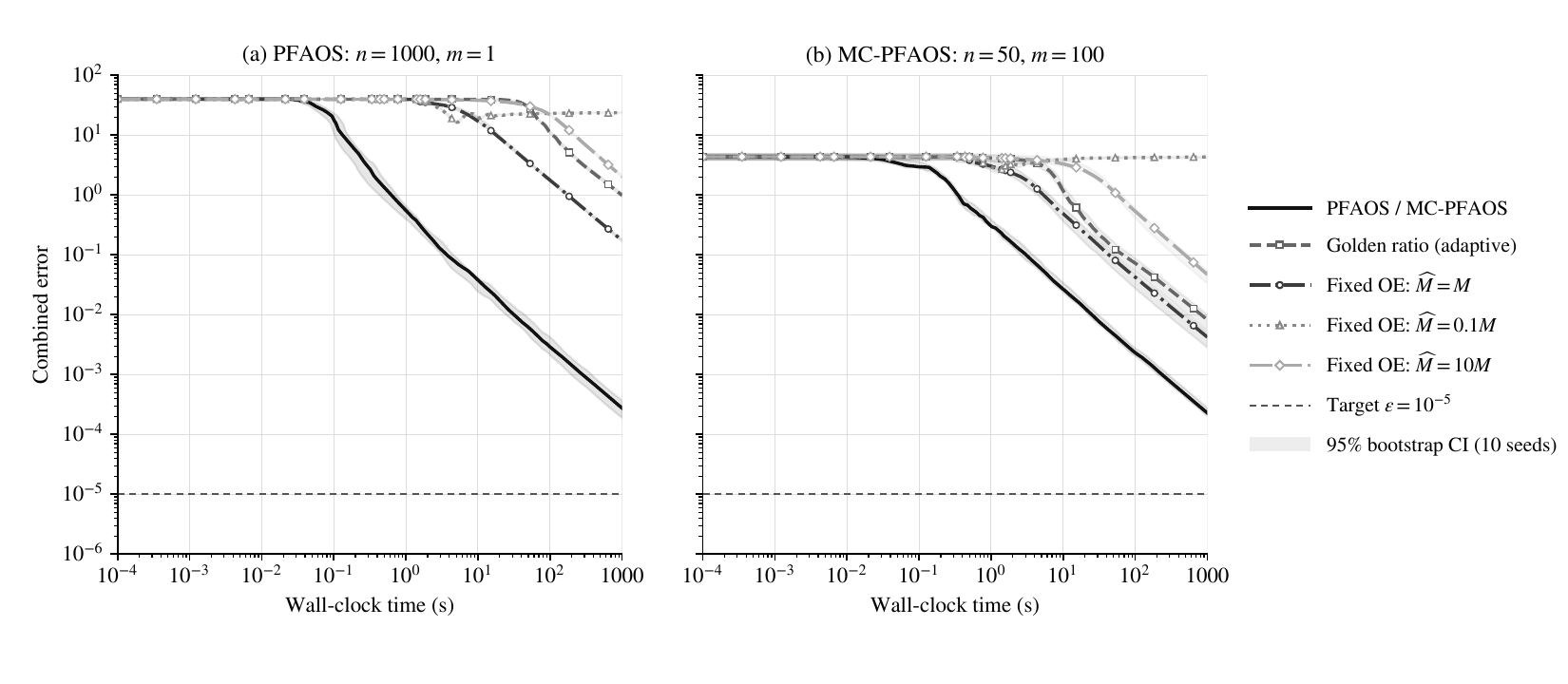}
  \caption{Combined error 
    versus solver time for (a) PFAOS and (b) MC-PFAOS,
    compared with OE and the Golden Ratio Algorithm.
    Here $\widehat{\operatorname{Gap}}$ denotes the numerical gap estimate.
    Central curves show medians over 10 random instances, and shaded bands show 95\% bootstrap confidence intervals for the medians. The dashed line marks the target accuracy $10^{-5}$.}
  \label{fig:pfaos_gap_n10000}
\end{figure}


 Figure \ref{fig:pfaos_gap_n10000} shows a substantial computational advantage for PFAOS and MC-PFAOS over the projection-based benchmarks.
 Each configuration uses 10 paired random seeds. 
 After \(1000\) seconds, PFAOS attains a median estimated weak gap in order of \( 10^{-4}\), compared with \( 10^{-1}\) for OE initialized at the true operator constant \(M\).
 MC-PFAOS achieves a zero estimated weak gap, although its remaining constraint violation yields a median combined error of \(10^{-4}\). 
 OE is sensitive to its fixed smoothness estimate: using \(10M\) leads to failed convergence substantially, whereas using \(0.1M\) leads to persistently large errors. 
 The adaptive golden-ratio benchmark also reduces the gap but converges considerably more slowly than the proposed methods. These comparisons reflect the computational savings of PFAOS from  avoiding full projection onto the nonlinear feasible set at every step.

\begin{table*}[!htbp]
\centering
\caption{First observed valid joint numerical-error threshold times (seconds) and terminal constraint violation (m=1) over 10 attempted runs per configuration, with a 1000-second method cap.}
\label{tab:single_constraint_comparative}
\setlength{\tabcolsep}{3pt}
\resizebox{\textwidth}{!}{
\begin{tabular}{clccccc}
\hline
$n$ & Method & $\widehat E\le10^{-2}$ & $\widehat E\le10^{-3}$ & $\widehat E\le10^{-4}$ & $\widehat E\le10^{-5}$ & Final $V$ \\
\hline
10 & PFAOS & \shortstack{$10/10$\\0.019/0.050/0.081} & \shortstack{$10/10$\\0.051/0.475/1.03} & \shortstack{$10/10$\\0.152/6.37/12.52} & \shortstack{$10/10$\\0.834/66.92/134.57} & 0.00e+00/9.81e-07/3.75e-06 \\
10 & OE ($\widehat M=M$) & \shortstack{$10/10$\\1.09/3.32/5.91} & \shortstack{$10/10$\\5.91/30.85/55.31} & \shortstack{$10/10$\\62.35/305.87/470.33} & \shortstack{$2/10$\\657.32/691.59/725.85} & 0.00e+00/0.00e+00/0.00e+00 \\
10 & Golden Ratio & \shortstack{$10/10$\\2.30/6.46/9.95} & \shortstack{$10/10$\\17.38/63.72/104.91} & \shortstack{$10/10$\\140.24/629.26/1000.00} & \shortstack{$0/10$\\--} & 0.00e+00/0.00e+00/0.00e+00 \\
\hline
50 & PFAOS & \shortstack{$10/10$\\0.039/0.167/0.345} & \shortstack{$10/10$\\0.658/1.68/3.03} & \shortstack{$10/10$\\5.88/16.80/33.00} & \shortstack{$10/10$\\65.57/167.61/273.30} & 9.65e-07/1.57e-06/3.13e-06 \\
50 & OE ($\widehat M=M$) & \shortstack{$10/10$\\4.97/8.97/13.53} & \shortstack{$10/10$\\37.51/91.32/142.63} & \shortstack{$9/10$\\395.51/793.32/1000.00} & \shortstack{$0/10$\\--} & 0.00e+00/0.00e+00/0.00e+00 \\
50 & Golden Ratio & \shortstack{$10/10$\\7.86/15.51/22.94} & \shortstack{$10/10$\\79.12/153.55/241.89} & \shortstack{$3/10$\\834.07/902.79/1000.00} & \shortstack{$0/10$\\--} & 0.00e+00/0.00e+00/0.00e+00 \\
\hline
100 & PFAOS & \shortstack{$10/10$\\0.134/0.243/0.373} & \shortstack{$10/10$\\1.53/2.44/3.27} & \shortstack{$10/10$\\17.43/24.35/34.03} & \shortstack{$10/10$\\148.46/237.12/299.46} & 1.54e-06/2.29e-06/4.20e-06 \\
100 & OE ($\widehat M=M$) & \shortstack{$10/10$\\21.54/32.03/43.67} & \shortstack{$10/10$\\227.20/300.77/392.63} & \shortstack{$0/10$\\--} & \shortstack{$0/10$\\--} & 0.00e+00/0.00e+00/0.00e+00 \\
100 & Golden Ratio & \shortstack{$10/10$\\39.49/53.72/64.73} & \shortstack{$10/10$\\416.46/566.40/682.52} & \shortstack{$0/10$\\--} & \shortstack{$0/10$\\--} & 0.00e+00/0.00e+00/0.00e+00 \\
\hline
200 & PFAOS & \shortstack{$10/10$\\0.315/0.475/0.607} & \shortstack{$10/10$\\3.07/5.03/7.28} & \shortstack{$10/10$\\34.54/52.09/73.01} & \shortstack{$10/10$\\323.91/505.63/661.39} & 1.98e-06/2.59e-06/3.34e-06 \\
200 & OE ($\widehat M=M$) & \shortstack{$10/10$\\137.31/184.11/206.46} & \shortstack{$0/10$\\--} & \shortstack{$0/10$\\--} & \shortstack{$0/10$\\--} & 0.00e+00/0.00e+00/0.00e+00 \\
200 & Golden Ratio & \shortstack{$10/10$\\259.07/308.81/367.15} & \shortstack{$0/10$\\--} & \shortstack{$0/10$\\--} & \shortstack{$0/10$\\--} & 0.00e+00/0.00e+00/0.00e+00 \\
\hline
500 & PFAOS & \shortstack{$10/10$\\0.556/0.784/1.22} & \shortstack{$10/10$\\5.88/8.08/10.79} & \shortstack{$10/10$\\65.13/81.18/122.52} & \shortstack{$10/10$\\554.19/804.24/998.11} & 3.16e-06/4.37e-06/7.05e-06 \\
500 & OE ($\widehat M=M$) & \shortstack{$0/10$\\--} & \shortstack{$0/10$\\--} & \shortstack{$0/10$\\--} & \shortstack{$0/10$\\--} & 0.00e+00/0.00e+00/0.00e+00 \\
500 & Golden Ratio & \shortstack{$0/10$\\--} & \shortstack{$0/10$\\--} & \shortstack{$0/10$\\--} & \shortstack{$0/10$\\--} & 0.00e+00/0.00e+00/0.00e+00 \\
\hline
1000 & PFAOS & \shortstack{$10/10$\\2.09/2.77/4.77} & \shortstack{$10/10$\\18.29/27.20/44.92} & \shortstack{$10/10$\\186.59/270.60/410.31} & \shortstack{$0/10$\\--} & 8.83e-06/1.33e-05/2.34e-05 \\
1000 & OE ($\widehat M=M$) & \shortstack{$0/10$\\--} & \shortstack{$0/10$\\--} & \shortstack{$0/10$\\--} & \shortstack{$0/10$\\--} & 0.00e+00/0.00e+00/0.00e+00 \\
1000 & Golden Ratio & \shortstack{$0/10$\\--} & \shortstack{$0/10$\\--} & \shortstack{$0/10$\\--} & \shortstack{$0/10$\\--} & 0.00e+00/0.00e+00/0.00e+00 \\
\hline
\end{tabular}}
\smallskip
\parbox{\textwidth}{\footnotesize\textit{Note:}   Each threshold cell shows observed valid successes $s/N$, then best/mean/worst times over those successful runs only.   -- means no observed valid crossing before the recorded termination. Final $V$ gives best/mean/worst over all terminal   outputs. Times are  in seconds.  }
\end{table*}

The single-constraint study fixes $m=1$ and uses
$n\in\{10,50,100,200,500,1000\}$, with ten paired random seeds per dimension.
For the experiments in Table~\ref{tab:single_constraint_comparative},
PFAOS reaches the joint numerical-error threshold $10^{-4}$ in all 60 runs
and reaches $10^{-5}$ in every run with $n\le500$.
At $n=1000$, all ten runs reach $10^{-4}$, but none reaches $10^{-5}$
within the 1000-second budget.
For thresholds attained by both PFAOS and a baseline, PFAOS has the smaller
conditional mean hitting time.
In contrast, fixed OE with $\widehat M=M$  reaches $10^{-5}$ in only two of the ten $n=10$ runs, and the Golden Ratio method also records no $10^{-5}$ crossing.
These results demonstrate PFAOS's stronger empirical performance regarding both the the VI gap and constraint violations.

\section{Conclusion}\label{sec:Conclusion}
This work develops PFAOS and MC-PFAOS for monotone VIs with convex smooth functional constraints.
Both methods separate operator evaluations from constraint-gradient
or Jacobian evaluations, while avoiding projections onto the complex
feasible set.
Neither method requires prior knowledge of the Lipschitz constant
of the operator or the smoothness constants of the constraints;
MC-PFAOS also requires no prior uniform bound on the Jacobian.
Under $\Lambda\geq\|\lambda^*\|+1$, both methods require at most
$\mathcal O(1/\sqrt{\varepsilon})$ constraint-gradient or Jacobian
evaluations and $\mathcal O(1/\varepsilon)$ operator evaluations,
including evaluations in rejected line-search trials, to obtain an $\varepsilon$-approximate solution.
With the special first-iteration search and automatic stopping
rule, both methods accommodate arbitrarily small positive initial
estimates and certify an $\varepsilon$-approximate
solution without requiring knowledge of the diameter of the set $X$, while preserving these
optimal complexity bounds.
\ACKNOWLEDGMENT{
This work is partially supported by the AFOSR grant FA9550-25-1-0278.
During the preparation of this work, all theoretical derivations and numerical algorithms were developed by the authors without the use of AI tools. The authors used ChatGPT to improve the readability of the manuscript. After using this tool, the authors reviewed and edited the content as needed and take full responsibility for the content of the publication. The authors also used ChatGPT to debug the implementation code for the numerical experiments and reviewed the code as needed to ensure that it faithfully implements the proposed algorithms. This work was supported in part by OpenAI API credits provided by Clemson University and administered by Clemson University Research Computing and Data (RCD).}

\bibliographystyle{informs2014-no-url-note}
\bibliography{yan, JoC_submission/sample} 

\newpage

\begin{APPENDICES}
\renewcommand{\theHsection}{appendix.\Alph{section}}
\section{Proofs for Section~\ref{sec:mc_pfaos}}
\label{app:mc_pfaos}

\begin{proof}{Proof of Proposition~\ref{thm:mc_pfaos_inner_analysis}.}
Fix any $t\ge1$. Applying Lemma~\ref{thm:three_point} to the primal and dual subproblems in \eqref{eq:mc_pfaos_subproblem}, we have 
\begin{equation*}
    \begin{aligned}
       0\le  &\langle\ell_k^s(x_k^{t,s}),\lambda_k^{t,s}-\lambda\rangle
    -\tfrac{\eta_k^s+p_k^{t,s}}{2}\|\lambda-\lambda_k^{t,s}\|^2
    +
    \tfrac{\eta_k^s}{2}
    \|\lambda-\lambda_{k-1}\|^2
    +\tfrac{p_k^{t,s}}{2}
    (\|\lambda-\lambda_k^{t-1,s}\|^2-\|\lambda_k^{t,s}-\lambda_k^{t-1,s}\|^2), \text{ and }\\
        0\le& -\langle
        F(x_k^{t-1,s})+\theta_k^{t,s}\Delta_k^{t-1,s}+A_k^{t,s},
         x_k^{t,s}-x
         \rangle
        -\tfrac{\eta_k^s+p_k^{t,s}}{2}\|x-x_k^{t,s}\|^2
        \\
        &+
        \tfrac{\eta_k^s}{2}
        (\|x-x_{k-1}\|^2-\|x_k^{t,s}-x_{k-1}\|^2)
        +\tfrac{p_k^{t,s}}{2}
        (\|x-x_k^{t-1,s}\|^2-\|x_k^{t,s}-x_k^{t-1,s}\|^2).
    \end{aligned}
\end{equation*}
Combining the above relations with  
\(
\langle F(x),x_k^{t,s}-x\rangle
\le
\langle F(x_k^{t,s}),x_k^{t,s}-x\rangle
\) from the monotonicity of $F$, and applying the Cauchy--Schwarz inequality and rearranging terms yields
\(
    \langle F(x), x_k^{t,s}-x\rangle 
    \le 
    -\langle A_k^{t,s},x_k^{t,s}-x\rangle + B_k^{t,s}+ C_k^{t,s} + D_k^{t,s},
\)
where 
\begin{equation*}
    \begin{aligned}
        B_k^{t,s} := & \langle \Delta_k^{t,s}, x_k^{t,s}-x\rangle - \theta_k^{t,s}\langle\Delta_k^{t-1,s}, x_k^{t-1,s}-x\rangle   + \tfrac{p_k^{t,s}}{2}(\|x-x_k^{t-1,s}\|^2+\|\lambda-\lambda_k^{t-1,s}\|^2)
        \\
        & - \tfrac{\eta_k^s+p_k^{t,s}}{2}(\|x-x_k^{t,s}\|^2+\|\lambda-\lambda_k^{t,s}\|^2)
        \\
        C_k^{t,s} : = & \theta_k^{t,s}\|\Delta_k^{t-1,s}\|\|x_k^{t,s}-x_k^{t-1,s}\| - \tfrac{p_k^{t,s}}{2}\|x_k^{t,s}-x_{k}^{t-1,s}\|^2
        \\
        D_k^{t,s} := & \tfrac{\eta_k^s}{2}(\|x-x_{k-1}\|^2 - \|x_k^{t,s}-x_{k-1}\|^2 + \|\lambda-\lambda_{k-1}\|^2) -\tfrac{p_k^{t,s}}{2}\|\lambda_k^{t,s}-\lambda_k^{t-1,s}\|^2
        +\langle \ell_k^s(x_k^{t,s}),\lambda_k^{t,s}-\lambda\rangle.
    \end{aligned}
\end{equation*}
Multiplying both sides of the above relation by $1/M_k^{t,s}$, and summing over $t=1,\ldots,T_k^s$, while considering the definition of $\xhat_k^s$ gives immediately that
\begin{equation}\label{eq:mcpfaos_inner_sum}
\begin{aligned}
    \mathcal{B}_k^s\langle F(x), \xhat_k^s-x\rangle 
    \le 
    -\sum\nolimits_{t=1}^{T_k^s}\tfrac{1}{M_k^{t,s}}\langle A_k^{t,s},x_k^{t,s}-x\rangle 
    + \sum\nolimits_{t=1}^{T_k^s} B_k^{t,s}/M_k^{t,s}
    + \sum\nolimits_{t=1}^{T_k^s}C_k^{t,s}/M_k^{t,s} 
    + \sum\nolimits_{t=1}^{T_k^s}D_k^{t,s}/M_k^{t,s}.
\end{aligned}
\end{equation}
Now, we bound the last three sums in the right hand side of \eqref{eq:mcpfaos_inner_sum} one by one. First, the parameter settings in Algorithm~\ref{alg:mc-pfaos} ensure that $\theta_k^{t,s} = M_k^{t,s}/M_k^{t-1,s}$ and $p_k^{t,s}/M_k^{t,s}=(p_k^{t-1,s}+\eta_k^s)/M_k^{t-1,s}$, $\forall t\ge2$. Thus, canceling out the recursive terms while considering the definitions of $M_k^{s}$, $x_k^s$, $x_k^{0,s}$, we have $\sum\nolimits_{t=1}^{T_k^s}B_k^{t,s}/M_k^{t,s}=B_k^s$. Second,  Algorithm \ref{alg:mc-pfaos} also indicates that $p_k^{t,s}/M_k^{t,s}\ge3$, $\forall t\ge1$, so
\begin{equation*}
\|x_k^{t-1,s}-x_k^{t-2,s}\|\|x_k^{t,s}-x_k^{t-1,s}\| -\tfrac{p_k^{t,s}}{6M_k^{t,s}}\|x_k^{t,s}-x_k^{t-1,s}\|^2 - \tfrac{p_k^{t-1,s}}{6M_
    k^{t-1,s}}\|x_k^{t-1,s}-x_k^{t-2,s}\|^2\le0, \quad \forall t\ge 2.
\end{equation*}
Thus, given the relations 
    \(\|\Delta_k^{t-1,s}\|\le M_k^{t-1,s}\|x_k^{t-1,s}-x_k^{t-2,s}\|
    \)
    and 
    $\theta_k^{t,s} = M_k^{t,s}/M_k^{t-1,s}$, $\forall t\ge2$ 
    guaranteed by the algorithm, and the definitions of $x_k^{0,s}$ and $x_k^s$, we have $\sum\nolimits_{t=1}^{T_k^s}C_k^{t,s}/M_k^{t,s} \le C_k^s$. Finally, with the help of the convexity of $\|\cdot\|^2$, the invariant relation $\mathcal{B}_k^s\eta_k^s=1$, together with the definitions of $\hat{x}_k^s$, it is easy to see that $\sum\nolimits_{t=1}^{T_k^s}D_k^{t,s}/M_k^{t,s} \le D_k^s$. \hfill \halmos
\end{proof}

\begin{proof}{Proof of Proposition~\ref{thm:mc_pfaos_middle_gap}.}
Immediately from the definitions of $\xub_k$, $\hat{x}_k$, $\overline\lambda_k$ and $\lamhat_k$, it follows that
\begin{equation}\label{eq:onestep_outer_mcpfaos_overall}
    \begin{aligned}
        \mathcal G_k(x,\lambda) -(1-\gamma_k^{S_k})\mathcal G_{k-1}(x,\lambda)
=
\gamma_k^{S_k}\langle F(x),\xhat_k-x\rangle
+\langle G(\xub_k),\lambda\rangle
-(1-\gamma_k^{S_k})\langle G(\xub_{k-1}),\lambda\rangle
-\gamma_k^{S_k}\langle G(x),\hat\lambda_k\rangle.
    \end{aligned}
\end{equation}
Here, we make two observations. 
First, the acceptance test described at line~\ref{line:mcpfaos_middle_test} gives
\(
\|G(\xub_k)-\ell_k^{S_k}(\xub_k)\|
\le
\tfrac{\eta_k^{S_k}}{2\Lambda\gamma_k^{S_k}}
\|\xub_k-\xlb_k^{S_k}\|^2.
\) Second, 
since $\ell_k^{S_k}$ is affine and each $g_i$ is convex, we have
\(
\ell_k^{S_k}(\xub_k)
=(1-\gamma_k^{S_k})\ell_k^{S_k}(\xub_{k-1})
+\gamma_k^{S_k}\ell_k^{S_k}(\xhat_k),
\)
while
\(
G(\xub_{k-1})\ge\ell_k^{S_k}(\xub_{k-1})
\)
and
\(
G(x)\ge\ell_k^{S_k}(x)
\)
componentwise. Hence, using $\lambda,\hat\lambda_k\ge0$, Cauchy-Schwarz inequality and the relation
\(
\xub_k-\xlb_k^{S_k}=\gamma_k^{S_k}(\xhat_k-x_{k-1})
\) ensured by the algorithm, we obtain
\begin{equation}\label{eq:onestep_outer_mcpfaos_dual_terms}
        \langle G(\xub_k),\lambda\rangle
-(1-\gamma_k^{S_k})\langle G(\xub_{k-1}),\lambda\rangle
-\gamma_k^{S_k}\langle G(x),\hat\lambda_k\rangle
\le
\gamma_k^{S_k}
(
\langle\ell_k^{S_k}(\xhat_k),\lambda\rangle
-\langle\ell_k^{S_k}(x),\hat\lambda_k\rangle
+\tfrac{\eta_k^{S_k}\|\lambda\|}{2\Lambda}
\|\xhat_k-x_{k-1}\|^2)
\end{equation}
To reach the conclusion, we also need the following five facts. First, we have
\begin{equation}\label{eq:mc_pfaos_affine_cancellation}
\begin{aligned}
&\tfrac{1}{\mathcal B_k^{S_k}}
\sum\nolimits_{t=1}^{T_k^{S_k}}\tfrac{1}{M_k^{t,S_k}}
 (
\langle\ell_k^{S_k}(x_k^{t,S_k}),\lambda_k^{t,S_k}-\lambda\rangle
-\langle J(\xlb_k^{S_k})(x_k^{t,S_k}-x),\lambda_k^{t,S_k}\rangle
 )\\
=&\tfrac{1}{\mathcal B_k^{S_k}}
\sum\nolimits_{t=1}^{T_k^{S_k}}\tfrac{1}{M_k^{t,S_k}}
 (
\langle\ell_k^{S_k}(x),\lambda_k^{t,S_k}\rangle
-\langle\ell_k^{S_k}(x_k^{t,S_k}),\lambda\rangle
 )=\langle\ell_k^{S_k}(x),\hat\lambda_k\rangle
-\langle\ell_k^{S_k}(\xhat_k),\lambda\rangle.
\end{aligned}
\end{equation}
Here the first equality holds due to
$\ell_k^{S_k}(x_k^{t,S_k})-J(\xlb_k^{S_k})(x_k^{t,S_k}-x)
=\ell_k^{S_k}(x)$, while the second one comes from the
definitions of $\xhat_k$ and $\hat\lambda_k$ and the affinity of
$\ell_k^{S_k}$.
Second, directly from the definitions of $\lamtil_k^{t,S_k}$ and the relation $\theta_k^{t,s} = M_k^{t,s}/M_k^{t-1,s}$, $\forall t\ge2$,  given by the algorithm, it follows that 
\begin{equation}\label{eq:mc_pfaos_dual_telescoping}
\begin{aligned}
&\sum\nolimits_{t=2}^{T_k^{S_k}} 
(
\tfrac{1}{M_k^{t,S_k}}
\langle J(\xlb_k^{S_k})(x_k^{t,S_k}-x),
\lambda_k^{t,S_k}-\lamtil_k^{t,S_k}\rangle
+
\tfrac{\theta_k^{t,S_k}}{M_k^{t,S_k}}
\langle J(\xlb_k^{S_k})(x_k^{t,S_k}-x_k^{t-1,S_k}),
\lambda_k^{t-1,S_k}-\lambda_k^{t-2,S_k}\rangle
 )\\
&=\tfrac{1}{M_k^{S_k}}
\langle J(\xlb_k^{S_k})(x_k-x),
\lambda_k-\lambda_k^{T_k^{S_k}-1,S_k}\rangle-
\tfrac{1}{M_k^{1,S_k}}
\langle J(\xlb_k^{S_k})(x_k^{1,S_k}-x),
\lambda_k^{1,S_k}-\lambda_{k-1}\rangle,
\end{aligned}
\end{equation}
where the last equality also uses the definitions of $M_k^{S_k}$, $x_k$, $\lambda_k$, and
$\lambda_k^{0,S_k}$. This identity also holds when
$T_k^{S_k}=1$, with the sums interpreted as zero.
Third, reindexing terms yields
\begin{equation}\label{eq:mc_pfaos_dual_reindex}
\begin{aligned}
&-\sum\nolimits_{t=1}^{T_k^{S_k}}
\tfrac{p_k^{t,S_k}}{2M_k^{t,S_k}}
\|\lambda_k^{t,S_k}-\lambda_k^{t-1,S_k}\|^2
=-\tfrac{p_k^{T_k^{S_k},S_k}}{2M_k^{S_k}}
\|\lambda_k-\lambda_k^{T_k^{S_k}-1,S_k}\|^2
-\sum\nolimits_{t=2}^{T_k^{S_k}}
\tfrac{p_k^{t-1,S_k}}{2M_k^{t-1,S_k}}
\|\lambda_k^{t-1,S_k}-\lambda_k^{t-2,S_k}\|^2.
\end{aligned}
\end{equation}
Fourth, by utilizing Cauchy--Schwarz inequality, the
Jacobian--dual line-search condition described at line \ref{line:mcpfaos_t_loop} and the relations $\theta_k^{t,S_k}/M_k^{t,S_k}=1/M_k^{t-1,S_k}$ and $p_k^{t,S_k}/M_k^{t,S_k}\ge3$, we have, for $t\ge2$, 
\begin{equation}
\begin{aligned}
&-\tfrac{\theta_k^{t,S_k}}{M_k^{t,S_k}}
\langle J(\xlb_k^{S_k})(x_k^{t,S_k}-x_k^{t-1,S_k}),
\lambda_k^{t-1,S_k}-\lambda_k^{t-2,S_k}\rangle
-
\tfrac{p_k^{t-1,S_k}}{2M_k^{t-1,S_k}}
\|\lambda_k^{t-1,S_k}-\lambda_k^{t-2,S_k}\|^2
-
\tfrac{p_k^{t,S_k}}{6M_k^{t,S_k}}
\|x_k^{t,S_k}-x_k^{t-1,S_k}\|^2
\\
&\le
\|\lambda_k^{t-1,S_k}-\lambda_k^{t-2,S_k}\|
\|x_k^{t,S_k}-x_k^{t-1,S_k}\|
-\tfrac32\|\lambda_k^{t-1,S_k}-\lambda_k^{t-2,S_k}\|^2
-\tfrac12\|x_k^{t,S_k}-x_k^{t-1,S_k}\|^2
\le0.
\end{aligned}
\end{equation}
Lastly, using  the Cauchy--Schwarz inequality gives
\begin{equation}\label{eq:mc_pfaos_first_dual_boundary}
\begin{aligned}
&-\langle J(\xlb_{k-1}^{S_{k-1}})(x_k^{1,S_k}-x),
\lambda_{k-1}-\lambda_k^{-1,S_k}\rangle\\
\le&-\langle J(\xlb_{k-1}^{S_{k-1}})(x_{k-1}-x),
\lambda_{k-1}-\lambda_k^{-1,S_k}\rangle
+
\|J(\xlb_{k-1}^{S_{k-1}})^\top
(\lambda_{k-1}-\lambda_k^{-1,S_k})\|
\|x_k^{1,S_k}-x_{k-1}\|.
\end{aligned}
\end{equation}

Substitute \eqref{eq:conclusion_mcpfaos_inner} from Proposition~\ref{thm:mc_pfaos_inner_analysis} at $s=S_k$ and
\eqref{eq:onestep_outer_mcpfaos_dual_terms} into  \eqref{eq:onestep_outer_mcpfaos_overall}.
Then \eqref{eq:mc_pfaos_outer_gap} follows from  
\eqref{eq:mc_pfaos_affine_cancellation}--\eqref{eq:mc_pfaos_first_dual_boundary}, $\mathcal B_k^{S_k}\eta_k^{S_k}=1$, the inner-loop initializations and the definitions of the accepted outer iterates. \hfill \halmos
\end{proof}

\begin{proof}{Proof of Theorem~\ref{thm:mc_pfaos_global_convergence}.}
Divide \eqref{eq:mc_pfaos_outer_gap} by $\Gamma_k^{S_k}$ and sum over
$k=1,\ldots,N$. Since $\gamma_1^{S_1}=1$ and
$\Gamma_k^{S_k}=\Gamma_{k-1}^{S_{k-1}}(1-\gamma_k^{S_k})$ for $k\geq2$,
the left-hand side telescopes:
\begin{equation}\label{eq:mc_pfaos_outer_gap_sum}
\tfrac {\mathcal G_N(x,\lambda)}{\Gamma_N^{S_N}}
=\sum\nolimits_{k=1}^N
\tfrac {\mathcal G_k(x,\lambda)
-(1-\gamma_k^{S_k})\mathcal G_{k-1}(x,\lambda)}{\Gamma_k^{S_k}}
\le
\sum\nolimits_{k=1}^N (
E_k 
+F_k
+R'_k(\lambda)).
\end{equation}

Here, the expression of $E_k$ in \eqref{eq:mc_pfaos_outer_gap} is structured so that its terms telescope
over the outer iterations. To see this, the definition of $p_k^{t,S_k}$ gives
\(
(\eta_k^{S_k}+p_k^{T_k^{S_k},S_k})/M_k^{S_k}
=3+\eta_k^{S_k}\sum\nolimits_{t=1}^{T_k^{S_k}}1/M_k^{t,S_k}
=3+\eta_k^{S_k}\mathcal{B}_k^{S_k}
=4
\)
and \(
p_k^{1,S_k}/M_k^{1,S_k}=3
\).
Moreover, for $k\ge2$, the parameter settings 
imply
\(
\theta_k^{1,S_k}/M_k^{1,S_k}
=1/M_{k-1}^{S_{k-1}}
\)
,
\(
\Delta_k^{0,S_k}
=\Delta_{k-1}^{T_{k-1}^{S_{k-1}},S_{k-1}}
\).
Using $\Delta_1^{0,S_1}=0$ and
$\lambda_1^{0,S_1}-\lambda_1^{-1,S_1}=0$ yields: 
\begin{equation}\label{eq:mc_pfaos_sum_E}
\begin{aligned}
\sum\nolimits_{k=1}^N E_k
\le{}&\tfrac {1}{M_N^{S_N}}
\langle\Delta_N^{T_N^{S_N},S_N},x_N-x\rangle
+
\tfrac {1}{M_N^{S_N}}
\langle J(\xlb_N^{S_N})(x_N-x),
\lambda_N-\lambda_N^{T_N^{S_N}-1,S_N}\rangle\\
&+2(\|x-x_0\|^2+\|\lambda-\lambda_0\|^2
-\|x-x_N\|^2
).
\end{aligned}
\end{equation}
Furthermore, the final inner-loop acceptance conditions in line \ref{line:mcpfaos_t_loop}
give
\(
\|\Delta_N^{T_N^{S_N},S_N}\|
\le M_N^{S_N}\|x_N-x_N^{T_N^{S_N}-1,S_N}\|
\) and 
\(
\|J(\xlb_N^{S_N})^\top
(\lambda_N-\lambda_N^{T_N^{S_N}-1,S_N})\|
\le M_N^{S_N}\|\lambda_N-\lambda_N^{T_N^{S_N}-1,S_N}\|
\).
Using
$p_N^{T_N^{S_N},S_N}/M_N^{S_N}\ge3$ and Cauchy--Schwarz inequality, we have, respectively,
\begin{equation*}
\begin{aligned}
&\tfrac {1}{M_N^{S_N}}
\langle\Delta_N^{T_N^{S_N},S_N},x_N-x\rangle
-\|x-x_N\|^2
\le 
\tfrac {p_N^{T_N^{S_N},S_N}}{6M_N^{S_N}}
\|x_N-x_N^{T_N^{S_N}-1,S_N}\|^2, 
\\&
\tfrac {1}{M_N^{S_N}}
\langle J(\xlb_N^{S_N})(x_N-x),
\lambda_N-\lambda_N^{T_N^{S_N}-1,S_N}\rangle-\|x-x_N\|^2
\le
\tfrac {p_N^{T_N^{S_N},S_N}}{2M_N^{S_N}}
\|\lambda_N-\lambda_N^{T_N^{S_N}-1,S_N}\|^2.
\end{aligned}
\end{equation*}
Thus, substituting the above two inequalities into \eqref{eq:mc_pfaos_sum_E} gives \begin{equation}\label{eq:mc_pfaos_sum_E_further}
\sum\nolimits_{k=1}^N E_k
\le
2(\|x-x_0\|^2+\|\lambda-\lambda_0\|^2)
+
\tfrac {p_N^{T_N^{S_N},S_N}}{6M_N^{S_N}}
\|x_N-x_N^{T_N^{S_N}-1,S_N}\|^2
+
\tfrac {p_N^{T_N^{S_N},S_N}}{2M_N^{S_N}}
\|\lambda_N-\lambda_N^{T_N^{S_N}-1,S_N}\|^2
.
\end{equation}

We now derive the
following facts to bound $\sum_{k=1}^N F_k$. First, by changing indices and using the
inner-loop initializations, we have
\begin{equation}\label{eq:mc_pfaos_reindex_primal}
\begin{aligned}
&\sum\nolimits_{k=1}^N
\tfrac {p_k^{T_k^{S_k},S_k}}{6M_k^{S_k}}
\|x_k-x_k^{T_k^{S_k}-1,S_k}\|^2
=
\tfrac {p_N^{T_N^{S_N},S_N}}{6M_N^{S_N}}
\|x_N-x_N^{T_N^{S_N}-1,S_N}\|^2
+
\sum_{k=2}^N
\tfrac {p_{k-1}^{T_{k-1}^{S_{k-1}},S_{k-1}}}{6M_{k-1}^{S_{k-1}}}
\|x_{k-1}-x_k^{-1,S_k}\|^2, 
\\
&\sum\nolimits_{k=1}^N
\tfrac {p_k^{T_k^{S_k},S_k}}{2M_k^{S_k}}
\|\lambda_k-\lambda_k^{T_k^{S_k}-1,S_k}\|^2
=\tfrac {p_N^{T_N^{S_N},S_N}}{2M_N^{S_N}}
\|\lambda_N-\lambda_N^{T_N^{S_N}-1,S_N}\|^2
+
\sum_{k=2}^N
\tfrac {p_{k-1}^{T_{k-1}^{S_{k-1}},S_{k-1}}}{2M_{k-1}^{S_{k-1}}}
\|\lambda_{k-1}-\lambda_k^{-1,S_k}\|^2.
\end{aligned}
\end{equation}
Second, for $k\ge2$, the final inner-loop acceptance conditions in line \eqref{line:mcpfaos_t_loop} imply
\(
\|\Delta_k^{0,S_k}\|
\le M_{k-1}^{S_{k-1}}
\|x_{k-1}-x_k^{-1,S_k}\|
\) 
and 
\(
\|J(\xlb_{k-1}^{S_{k-1}})^\top
(\lambda_{k-1}-\lambda_k^{-1,S_k})\|
\le M_{k-1}^{S_{k-1}}
\|\lambda_{k-1}-\lambda_k^{-1,S_k}\|.
\)
Thus, utilizing the relations
$\theta_k^{1,S_k}/M_k^{1,S_k}=1/M_{k-1}^{S_{k-1}}$,
$p_k^{1,S_k}/M_k^{1,S_k}=3$, and
$p_{k-1}^{T_{k-1}^{S_{k-1}},S_{k-1}}/M_{k-1}^{S_{k-1}}\ge3$ guaranteed by the Algorithm \ref{alg:mc-pfaos}, 
we obtain, for $k\ge2$,
\begin{equation}\label{eq:mc_pfaos_outer_cross_F}
\begin{aligned}
&\tfrac {\theta_k^{1,S_k}}{M_k^{1,S_k}}
\|\Delta_k^{0,S_k}\|\|x_k^{1,S_k}-x_{k-1}\|
-\tfrac {p_k^{1,S_k}}{6M_k^{1,S_k}}
\|x_k^{1,S_k}-x_{k-1}\|^2
-
\tfrac {p_{k-1}^{T_{k-1}^{S_{k-1}},S_{k-1}}}{6M_{k-1}^{S_{k-1}}}
\|x_{k-1}-x_k^{-1,S_k}\|^2
\le 0, 
\\&
\tfrac {\theta_k^{1,S_k}}{M_k^{1,S_k}}
\|J(\xlb_{k-1}^{S_{k-1}})^\top
(\lambda_{k-1}-\lambda_k^{-1,S_k})\|
\|x_k^{1,S_k}-x_{k-1}\|
-\tfrac {p_k^{1,S_k}}{6M_k^{1,S_k}}
\|x_k^{1,S_k}-x_{k-1}\|^2
-\tfrac {p_{k-1}^{T_{k-1}^{S_{k-1}},S_{k-1}}}{2M_{k-1}^{S_{k-1}}}
\|\lambda_{k-1}-\lambda_k^{-1,S_k}\|^2
\le0.
\end{aligned}
\end{equation}
Applying \eqref{eq:mc_pfaos_reindex_primal} and \eqref{eq:mc_pfaos_outer_cross_F},
using $\Delta_1^{0,S_1}=0$ and $\lambda_0-\lambda_1^{-1,S_1}=0$,
and dropping the nonpositive term
$-p_1^{1,S_1}\|x_1^{1,S_1}-x_0\|^2/(3M_1^{1,S_1})$, we obtain
\begin{equation*}\label{eq:mc_pfaos_sum_F}
\sum\nolimits_{k=1}^N F_k
\le-\tfrac {p_N^{T_N^{S_N},S_N}}{6M_N^{S_N}}
\|x_N-x_N^{T_N^{S_N}-1,S_N}\|^2
-\tfrac {p_N^{T_N^{S_N},S_N}}{2M_N^{S_N}}
\|\lambda_N-\lambda_N^{T_N^{S_N}-1,S_N}\|^2. \hfill \halmos
\end{equation*}
\end{proof}

\begin{proof}{Proof of Lemma~\ref{thm:mc_pfaos_Gamma_rate}.}
For any $u,v\in X$, the smoothness of the component functions gives
\begin{equation*}
\|G(v)-G(u)-J(u)(v-u)\|^2 
=
\sum\nolimits_{i=1}^m
[
g_i(v)-g_i(u)-\langle\nabla g_i(u),v-u\rangle
]^2
\leq
\tfrac{L^2}{4}\|v-u\|^4.
\end{equation*}
Hence, the middle-loop acceptance condition of MC-PFAOS is guaranteed
whenever
\(
\eta_k^s/\gamma_k^s\geq\Lambda L.
\)
Since the updates of $\eta_k^s$ and $\gamma_k^s$, and the definition of
$\Gamma_k^s$, are the same as those in PFAOS, the remaining proof of
Lemma~\ref{thm:pfaos_Gamma_rate} applies. \hfill \halmos
\end{proof}

\begin{proof}{Proof of Corollary~\ref{thm:mc_pfaos_duality_gap_rate}.}
For any $x\in\widetilde X$, since $G(x)\leq0$ and
$\overline\lambda_N\geq0$ componentwise,
Theorem~\ref{thm:mc_pfaos_global_convergence} and $\lambda_0=0$ give
\(
\langle F(x),\xub_N-x\rangle
\leq\mathcal G_N(x,0)
\leq2\Gamma_N^{S_N}D_X^2.
\)
Taking the supremum over $x\in\widetilde X$ gives the weak VI
gap bound in terms of $\Gamma_N^{S_N}$.
For the constraint violation, the claim is immediate if
$[G(\xub_N)]_+=0$. Otherwise, define
\(
u_N:=[G(\xub_N)]_+/\|[G(\xub_N)]_+\|.
\)
By \eqref{eq:standing_primal_dual_pair} and convexity of the
component functions,
\begin{equation*}
\langle F(x^*),\xub_N-x^*\rangle
\geq-\langle J(x^*)^\top\lambda^*,\xub_N-x^*\rangle
\geq-\langle\lambda^*,G(\xub_N)-G(x^*)\rangle
=-\langle\lambda^*,G(\xub_N)\rangle
\geq-\|\lambda^*\|\,\|[G(\xub_N)]_+\|.
\end{equation*}
Since $\Lambda u_N\in\mathbb R_+^m$,
$\|\Lambda u_N\|=\Lambda$, and
$\langle G(\xub_N),u_N\rangle=\|[G(\xub_N)]_+\|$,
using $G(x^*)\leq0$ and $\overline\lambda_N\geq0$ gives
\begin{equation*}
(\Lambda-\|\lambda^*\|)\|[G(\xub_N)]_+\|
\leq
\langle F(x^*),\xub_N-x^*\rangle
+\Lambda\langle G(\xub_N),u_N\rangle
\leq\mathcal G_N(x^*,\Lambda u_N)
\leq2\Gamma_N^{S_N}(D_X^2+\Lambda^2).
\end{equation*}
The explicit rates follow from
Lemma~\ref{thm:mc_pfaos_Gamma_rate}. \hfill \Halmos
\end{proof}

\begin{proof}{Proof of Proposition~\ref{thm:mc_pfaos_oracle_complexity}.}
Let $N_\varepsilon\geq1$ be the first outer iteration satisfying
\(
2\Gamma_{N_\varepsilon}^{S_{N_\varepsilon}}
(D_X^2+\Lambda^2)\leq\varepsilon.
\)
By Corollary~\ref{thm:mc_pfaos_duality_gap_rate},
$\xub_{N_\varepsilon}$ is an $\varepsilon$-approximate solution
of \eqref{eq:problem_of_interest} and
\(
N_\varepsilon
\leq1+\sqrt{\tfrac{32\Lambda L(D_X^2+\Lambda^2)}{\varepsilon}}.
\)
Each middle trial requires one evaluation of $J$.
Thus, the counting argument in
Proposition~\ref{thm:pfaos_oracle_complexity},  gives at most
\(
1+2N_\varepsilon+\lceil\log_2\tfrac{2L}{L_0}\rceil
\)
Jacobian evaluations, including initialization.
Using the bound on $N_\varepsilon$ yields the stated complexity.
For the operator evaluations, $\max_{x\in X}\|J(x)\|<\infty$
because $X$ is compact and $J$ is continuous.
Both inner line-search conditions hold whenever
$M_k^{t,s}\geq\max\{M,\max_{x\in X}\|J(x)\|\}$.
Thus, whenever backtracking occurs, the minimality of
$q_k^{t,s}$ gives
$M_k^{t,s}<2\max\{M,\max_{x\in X}\|J(x)\|\}$.
Moreover, since $M_0\leq2\Lambda L$ and
$\eta_k^s\leq\eta_k^{S_k}\le \eta_k^{S_k}/\gamma_k^{S_k}\leq4\Lambda L$,
the initialization
$M_k^{0,s}=\max\{\eta_k^s,M_k^{s-1}\}$
preserves the bound $\overline M_{\rm MC}$, while an inner
update without backtracking leaves the estimate unchanged.
Induction therefore yields $M_k^{t,s}\leq\overline M_{\rm MC}$
for every executed inner iteration.
Using this bound and the invariant $\eta_k^s\mathcal B_k^s=1$,
we repeat the counting argument in
Proposition~\ref{thm:pfaos_oracle_complexity}, now with the
stopping condition
\(
2\Gamma_{N_\varepsilon}^{S_{N_\varepsilon}}
(D_X^2+\Lambda^2)\leq\varepsilon,
\)
to obtain the stated operator complexity. \hfill \Halmos
\end{proof}

\section{Proofs for Section~\ref{sec:stopping_criteria}}
\label{app:stopping_criteria}

\begin{proof}{Proof of Proposition~\ref{thm:stopping_certificate}.}
Using relations $\Gamma_k^{S_k}=\gamma_k^{S_k}\eta_k^{S_k}$, $\eta_k^{S_k}\mathcal B_k^{S_k}=1$, $\gamma_1^{S_1}=1$ and $\Gamma_k^{S_k}=(1-\gamma_k^{S_k})\Gamma_{k-1}^{S_{k-1}}$, $k\geq2$ give
\begin{equation*}
\Gamma_N^{S_N}
\sum\nolimits_{k=1}^N\sum\nolimits_{t=1}^{T_k^{S_k}}\tfrac{1}{M_k^{t,S_k}}=\Gamma_N^{S_N}
\sum\nolimits_{k=1}^N\tfrac{1}{\eta_k^{S_k}}=1,\text{ and }
\Gamma_N^{S_N}
\sum\nolimits_{k=1}^N\sum\nolimits_{t=1}^{T_k^{S_k}}
\tfrac{x_k^{t,S_k}}{M_k^{t,S_k}}
=\Gamma_N^{S_N}
\sum\nolimits_{k=1}^N\mathcal{B}_k^{S_k}\xhat_k
=\xub_N.
\end{equation*}
Then, for any $x\in\widetilde X$, monotonicity of $F$ and convexity of
the constraints give
\begin{equation*}
\langle F(x),\xub_N-x\rangle
=\Gamma_N^{S_N}
\sum\nolimits_{k=1}^N\sum\nolimits_{t=1}^{T_k^{S_k}}
\tfrac{\langle F(x),x_k^{t,S_k}-x\rangle}{M_k^{t,S_k}}
\leq\Gamma_N^{S_N}
\sum\nolimits_{k=1}^N\sum\nolimits_{t=1}^{T_k^{S_k}}
\tfrac{\langle F(x_k^{t,S_k}),x_k^{t,S_k}-x\rangle}{M_k^{t,S_k}}
\leq\Phi_N(x),
\end{equation*}
where the last inequality uses
$\ell_k^{S_k}(x)\leq G(x)\leq0$ and
$\lambda_k^{t,S_k}\geq0$ componentwise.
Taking the supremum over $\widetilde X\subseteq X$ proves the
weak VI gap bound.

To bound the certificate, retain
$\langle F(x_k^{t,S_k}),x_k^{t,S_k}-x\rangle$ in the inner-loop
optimality inequalities, before applying monotonicity, and set the
test multiplier to $\lambda=0$.
Specifically, for PFAOS, the proof of
Proposition~\ref{thm:pfaos_inner_analysis}, with the affine constraint
term moved to the left, gives
\begin{equation*}
\sum\nolimits_{t=1}^{T_k^{S_k}}
\tfrac{
\langle F(x_k^{t,S_k}),x_k^{t,S_k}-x\rangle
-\lambda_k^{t,S_k}\ell_k^{S_k}(x)
}{M_k^{t,S_k}}
\leq
\tfrac{H_k+I_k}{\Gamma_k^{S_k}}
-\tfrac12\|\xhat_k-x_{k-1}\|^2,
\end{equation*}
where $H_k$ and $I_k$ are defined in
Proposition~\ref{thm:pfaos_middle_gap}.
Summing over $k$ and applying the telescoping and endpoint estimates
in the proof of Theorem~\ref{thm:pfaos_global_convergence} yields
\begin{equation*}
\tfrac{\Phi_N(x)}{\Gamma_N^{S_N}}
\leq
\sum\nolimits_{k=1}^N\tfrac{H_k+I_k}{\Gamma_k^{S_k}}
-\tfrac12\sum\nolimits_{k=1}^N\|\xhat_k-x_{k-1}\|^2
\leq\tfrac32\|x-x_0\|^2.
\end{equation*}
While for MC-PFAOS, retain the same operator term in
Proposition~\ref{thm:mc_pfaos_inner_analysis} and use the affine
cancellation in \eqref{eq:mc_pfaos_affine_cancellation}.
The resulting sum is bounded by
\begin{equation*}
\tfrac{\Phi_N(x)}{\Gamma_N^{S_N}}
\leq
 \sum\nolimits_{k=1}^N(E_k+F_k)  |_{\lambda=0}
-\tfrac12\sum\nolimits_{k=1}^N\|\xhat_k-x_{k-1}\|^2
\leq2\|x-x_0\|^2,
\end{equation*}
where the last inequality follows from the bounds on
$\sum_k E_k$ and $\sum_k F_k$ in the proof of
Theorem~\ref{thm:mc_pfaos_global_convergence} and $\lambda_0=0$.
Maximizing over $X$ and taking the positive part proves the stated
certificate bounds. The same arguments apply to retained trials of
Algorithm~\ref{alg:special_k1_pfaos} and \ref{alg:special_k1_mcpfaos}, as in
Lemma~\ref{thm:special_k1_convergence}.\hfill \Halmos
\end{proof}

\begin{proof}{Proof of Proposition~\ref{thm:stopping_complexity}.}
The accuracy claim follows directly from
Proposition~\ref{thm:stopping_certificate}.
Suppose now that $\Lambda\geq\|\lambda^*\|+1$.
By Proposition~\ref{thm:stopping_certificate} and
Corollaries~\ref{thm:pfaos_duality_gap_rate} and
\ref{thm:mc_pfaos_duality_gap_rate}, both conditions in
\eqref{eq:computable_stopping_rule} hold whenever
$\tfrac32\Gamma_N^{S_N}D_X^2\leq\varepsilon$ for PFAOS or
$2\Gamma_N^{S_N}(D_X^2+\Lambda^2)\leq\varepsilon$ for MC-PFAOS.
Thus, the rule stops no later than the thresholds used in the
original counting arguments.

The same implication holds for every trial at $k=1$ that passes the middle test.
If a trial fails the new stopping test, it also fails the original
accuracy test. Consequently, for a final tested index $s\geq2$,
the lower bounds
$\eta_1^s>\varepsilon/(3D_X^2)$ for PFAOS and
$\eta_1^s>\varepsilon/[4(D_X^2+\Lambda^2)]$ for MC-PFAOS
used in Lemma~\ref{thm:special_k1_cost} remain valid.
If the search continues to $k=2$, it still retains the trial
preceding a failed middle test, so $\eta_1^{S_1}<2\Lambda L$.
Hence, the first-iteration and subsequent counting arguments apply
unchanged. Certificate construction uses stored oracle values, and
$G(\xub_N)$ is already evaluated by the middle test, so the rule
requires no additional evaluations of $F$, $\nabla g$, or $J$. \hfill \Halmos
\end{proof}
\end{APPENDICES}

\end{document}